\documentclass[ijoc,nonblindrev]{oo} 

\OneAndAHalfSpacedXII 

\usepackage{natbib}
 \bibpunct[, ]{(}{)}{,}{a}{}{,}%
 \def\bibfont{\small}%

 \usepackage{float}
 \usepackage{enumitem}
 \usepackage{booktabs}
 \usepackage{multirow}
 \usepackage{subfigure}
 \usepackage[toc,page]{appendix}
\usepackage[utf8]{inputenc}
\usepackage{pgfplots}
\usepgfplotslibrary{groupplots,dateplot}
\usetikzlibrary{patterns,shapes.arrows}
\pgfplotsset{compat=newest}

 \usepackage{xcolor}

\TheoremsNumberedThrough     

\EquationsNumberedThrough    

\newcolumntype{L}[1]{>{\raggedright\let\newline\\\arraybackslash\hspace{0pt}}m{#1}}
\newcolumntype{C}[1]{>{\centering\let\newline\\\arraybackslash\hspace{0pt}}m{#1}}
\newcolumntype{R}[1]{>{\raggedleft\let\newline\\\arraybackslash\hspace{0pt}}m{#1}}                 
                 
\newcommand{\Graph}{G}
\newcommand{\nodeSet}{\mathcal{V}}
\newcommand{\node}{v}
\newcommand{\edgeSet}{\mathcal{L}}
\newcommand{\edge}{\ell}
\newcommand{\waveSet}{\Omega}
\newcommand{\wavelength}{\omega}
\newcommand{\LPSet}{\Lambda}
\newcommand{\lightpath}{\lambda}
\newcommand{\LPpath}{p}

\newcommand{\xirand}{\boldsymbol{\xi}}
\newcommand{\xival}{\xi}
\newcommand{\xiset}{\Xi}
\newcommand{\xisetSAA}{\hat{\Xi}}
\def \Dcurr {\bar{r}^{(s,d)}}
\def \Dnew {r^{(s,d)}(\xirand)}
\def \DnewSAA {r^{(s,d)}(\xival)}

\def \SD {\mathcal{SD}}
\newcommand{\SDcurr}[1]{\SD^1_{\textsc{#1}}}
\def \dual {\pi}
\def \SDsce {\mathcal{SD}^2(\xirand)}
\def \SDsceSAA {\mathcal{SD}^2(\xival)}

\def \nodefunc {\delta}
\def \link {x_{\edge,\wavelength}^{(s,d)}}
\def \linkhat {\hat{x}_{\edge,\wavelength}^{(s,d)}}
\def \linksce {y_{\edge,\wavelength}^{(s,d)}}
\def \linksceSAA {y_{\edge,\wavelength}^{(s,d)}}
\def \numgrantedsce {z^{(s,d)}}
\def \numgrantedsceSAA {z^{(s,d)}}
\def \numgrantedsceSAAhat {\hat{z}^{(s,d)}}
\newcommand {\Feas}[1] {\mathcal{F}^{\textsc{rwa}}(#1)}
\newcommand {\FeasR}[1] {\mathcal{F}^{\textsc{rwa}}_{\textsc{r}}(#1)}

\def \Z {\mathbb{Z}}
\def \P {\mathbb{P}}
\def \E {\mathbb{E}}
\def \N {\mathbb{N}}

\def \Q {\mathcal{Q}}
\def \cuts {\Gamma(\xival)}
\def \betacuts {\Upsilon(\xival)}
\def \dual {\pi}
\def \dualhat {\hat{\dual}}
\def \master {\textsc{master}}
\def \sub {\textsc{sub}}
\def \maxRWA {maxRWA}
\def \minRWA {minRWA}
\def \srwa {SmaxRWA}
\def \lr {SmaxLR}
\def \abilene {ABILENE}
\def \cost {COST239}
\def \nsf {NSF}
\def \atlanta {ATLANTA}
\def \usa {USA}
\def \brazil {BRAZIL}
\def \arrivals {\lambda^R}
\def \drops {\lambda^L}
\newcommand{\edgeSetw}{\mathcal{L}(\wavelength)}

\def \IPLP {IP-LP}
\def \IPIP {IP-IP}
\def \LPLP {LP-LP}

\newcommand{\BI}{\begin{itemize}}
\newcommand{\EI}{\end{itemize}}
\newcommand{\BE}{\begin{enumerate}}
\newcommand{\EE}{\end{enumerate}}

\newcommand{\BSE}{\begin{subequations}}
\newcommand{\ESE}{\end{subequations}}

\def \App {Appendix}

\newcommand{\revacc}[1]{\noindent{#1}}
\newcommand{\revminor}[1]{\noindent{#1}}
\newcommand{\rev}[1]{\noindent{#1}}

\usepackage{xr-hyper}
 \usepackage[colorlinks=true, linkcolor=blue, citecolor=blue]{hyperref}

\begin{document}


\RUNAUTHOR{Daryalal and Bodur}

\RUNTITLE{Stochastic RWA and Lightpath Rerouting in WDM Networks}

\TITLE{Stochastic RWA and Lightpath Rerouting in WDM Networks}

\ARTICLEAUTHORS{%
\AUTHOR{Maryam Daryalal, Merve Bodur}
\AFF{Department of Mechanical and Industrial Engineering, University of Toronto \\ \EMAIL{m.daryalal@mail.utoronto.ca}, \EMAIL{bodur@mie.utoronto.ca}
}
}

\ABSTRACT{%
In a telecommunication network, Routing and Wavelength Assignment (RWA) is the problem of finding lightpaths for incoming connection requests. When facing a dynamic traffic, greedy assignment of lightpaths to incoming requests based on predefined deterministic policies leads to a fragmented network that cannot make use of its full capacity due to stranded bandwidth. At this point service providers try to recover the capacity via a defragmentation process. We study this setting from two perspectives: ($ i $) while granting the connection requests via the RWA problem and ($ ii $) during the defragmentation process by lightpath rerouting. For both problems, we present the first two-stage stochastic integer programming model incorporating incoming request uncertainty to maximize the expected
grade of service. We develop a decomposition-based solution approach, which uses various relaxations of the problem and a newly developed problem-specific cut family.    
Simulation of two-stage policies for a variety of instances in a rolling-horizon framework of 52 stages shows that our stochastic models provide high-quality solutions compared to traditionally used deterministic ones. Specifically, the proposed provisioning policies yield improvements of up to 19\% in overall grade of service and 20\% in spectrum saving, while the stochastic lightpath rerouting policies grant up to 36\% more requests using up to just 4\% more bandwidth spectrum.

}%


\KEYWORDS{Two-stage Stochastic Programming, Routing and Wavelength Assignment, Lightpath Rerouting, Optical Networks}

\maketitle

\section{Introduction}
\label{sec:intro} 
The global internet traffic is steadily growing both in the number of users and the amount of transmitted data: it is predicted that by 2023, the total number of connected devices can reach as high as three times the world population \citep{cisco} and the amount of Internet traffic is multiplied by a factor of one thousand every twenty years \citep{wald2018impending}. This ever increasing traffic is in constant combat with the capacity of the information carriers, specifically in the context of \emph{optical networks}, the backbone of today's telecommunications systems \citep{majumdar2018optical}. While the cable operators have so far enjoyed the vast capacity provided by upgrading their existing submarine cables over the years, the industry is approaching the theoretical capacity of the optical fibers  \citep{cienaShannon2020}, called the Shannon limit \citep{shannon1948mathematical,essiambre2008capacity}. Combined with the 
increase in the number of users and the high bit-rate traffic requests such as online streaming services, this leads to the predicted \emph{capacity crunch} problem in telecommunication networks \citep{ellis2013are,wald2018impending, jara2020much, zhou2020combining}. Although capacity expansion by means of deploying fiber cables and extra equipment is a possible solution, ample investment and considerable amount of time required for the network expansion make it a strategic decision \citep{keiser1999review}.  The second and short-term option is efficient use of  existing capacity in optical networks that can postpone reaching the Shannon limit. Traditional studies aiming at enhancing this usage for a deployed network often assume a deterministic setting, which is in stark contrast with practice. In this paper, we illustrate the importance of incorporating traffic uncertainty into decision making, and introduce the first two-stage stochastic integer programming (2SP) model for the classical resource allocation problem in optical networks.

As the dominant transmission system in the telecommunications industry, \emph{wavelength division multiplexing} (WDM) networks are the subject of numerous studies aiming at improving the efficiency of optical networks. Therein, as the name of optical networks suggests, the transmission medium is an optical fiber carrying light. The key to their success is the ability to transmit multiple signals on different wavelengths through a single optical fiber, i.e., a WDM technology multiplexes and demultiplexes (joins and \rev{separates}) these signals at the nodes along a connection. 
In WDM networks, the pair of an optical fiber link and a wavelength is called a \emph{wavelink}. A connection is then set up via a \emph{lightpath}, which is defined as a specific frequency occupied on a sequence of optical fibers, i.e., a path  in the network together with a wavelength. 
The \emph{routing and wavelength assignment} (RWA)  problem is the provisioning problem in optical networks. Given a deployed WDM network and a set of connection requests as inputs, the RWA problem seeks to optimally assign one lightpath to each request considering two main conditions: $(i)$ every lightpath needs to use the same wavelength throughout its path, and $ (ii) $ two lightpaths with the same wavelength cannot share a link. WDM networks make it possible to efficiently scale the network \citep{cienawdm2020}, thus the RWA problem plays a vital role in the telecommunications industry that relies mostly on such transmission systems.

The RWA problem is quite challenging with a rich literature that has mainly focused on heuristic solutions. Exact solution approaches that can live up to the scale of demand the industry is facing today have only recently been proposed  \citep{jaumard2017efficient}. The objective function of the RWA problem can be defined based on different perspectives \citep{krishnaswamy2001algorithms}. One might consider minimizing the  resource usage while serving a set of connections that have to be granted, which corresponds to minimization of the number of used wavelinks (the \minRWA{} problem). Another possible objective function is to maximize the throughput measured as the number of granted requests  (the \maxRWA{} problem), which is equivalent to maximization of the grade of service (GoS), i.e., minimizing the blocking rate. For a brief overview of WDM networks, the RWA problem and the related literature, we refer the readers to \cite{kuri2003optimization} and \cite{daryalal2016efficient}.

Although in many studies on the RWA problem, traffic request is assumed to be static (and given in advance),  practical assumption in telecommunications industry, e.g. for a service provider, is a dynamic uncertain traffic. Mainly due to the immense size of  traffic requests in optical networks, stochasticity of the future demand has not been well leveraged at operational levels for deployed networks. Existing  works that consider stochasticity of traffic are mostly in the area of capacity planning and upgrade. \cite{kennington2003robust} develop a robust optimization model to determine the number of necessary equipment for a dense WDM network considering the future uncertain demand. For a survivable WDM network, \cite{leung2005capacity} propose a 2SP model that minimizes its deployment cost and possible augmentations in the future. More closely, \cite{lodha2007stochastic} study the light-trail design problem and formulate it as a 2SP model to decide on the set of light-trails to be configured. Light-trails are a different technology in WDM networks, where between every pair of nodes supporting the light-trail architecture, connections can share a wavelength. 
Considering a scenario tree with three scenarios and with the goal of minimizing the pre-planning costs,  \cite{kronberger2011impact} propose a three-stage stochastic programming model for dealing with the cases where future demand is strongly overestimated in the first stage. To determine the capital expenditure needed for capacity upgrade of a network, \cite{aparicio2012robust} present a robust optimization approach that considers stochastic demand with normal distribution. 

In the context of lightpath assignment and the RWA problem, although the stochasticity of traffic has been a discussion topic since the early studies \citep{ozdaglar2003routing},  explicitly incorporating  the uncertainty of the future traffic in granting decisions has been mostly avoided. Instead, the focus has been on the design of quasi-static/heuristic methods that tend to result in fewer congested links in the future.  To address this shortcoming, we propose a stochastic \maxRWA{} problem that, while assigning lightpaths, takes the future incoming traffic into account as a value function. To the best of our knowledge, this is the first study that proposes a stochastic model for lightpath assignment in a WDM network. Through numerical experiments, we  show that the use of the deterministic \maxRWA{} problem in a stochastic setting results in sub-optimal solutions  in terms of the number of granted requests. In all studied instances designed based on standard backbone networks,  the stochastic model can admit considerably more requests, up to 31 more connections, a 19\% improvement over the deterministic model. We note that as a result of traffic grooming and SONET/SDH technology, a single wavelength can carry multiple low bit-rate traffic connections with transmission rates as low as 51.84Mbps \citep{zhu2003review}; and today's metro and long-haul networks commonly operate at 100Gbps per wavelength \citep{INNISS201793}, and very recently at a remarkable bandwidth of 800Gbps \citep{ethernetalliance2020, cienaWave52020}.  
Therefore, even one more accepted request (assigned lightpath) is quite significant as it is equivalent to granting thousands of low bit-rate connections.

In a dynamic uncertain setting that traffic arrives and leaves in multiple time periods, a service provider sets up and tears down the lightpaths as needed. This type of traffic leads to network fragmentation \citep{ji2014dynamic}, meaning that there exists capacity stranded over the network, but it cannot be retrieved to build a lightpath for certain requests, resulting in high blocking rates. After it is identified that a network is highly fragmented, a \emph{network defragmentation} event is triggered. The state of the art in network defragmentation is to recover the stranded bandwidth by minimizing resource requirements modeled as a deterministic \minRWA{} problem  \citep{jaumard2019wavelength}. It obtains a set of lightpaths assigned to the granted connections on the network (target provisioning). A defragmentation process then involves rerouting of the lightpaths in a seamless migration fashion, i.e., in a way that a rerouted lightpath is set up before the original one is torn down. 
Noting that \minRWA{} ignores the future behavior of the traffic, as the second objective of this study, we question the effectiveness of \minRWA{} solution as the target provisioning in a defragmentation process. It is presumed that a target with fewer used resources performs better in the future for the incoming requests since available wavelinks are maximized \citep{jaumard2017efficient}.  In this paper, we propose to incorporate the natural uncertainty of incoming traffic in the design of the defragmentation target. We formulate the problem of finding this provisioning as a 2SP model and show its effectiveness in improving the performance of the network after multiple stages.

\noindent\textbf{Contributions}. Main contributions of this paper are summarized as follows:
\begin{itemize}
	\item We introduce the \emph{stochastic maxRWA} (\srwa{}) problem, a stochastic variant of  the \maxRWA{} problem incorporating traffic uncertainty, and formulate it as a 2SP model.
	\item Rather than using \minRWA{} as the target provisioning of defragmentation problem, we propose the \emph{stochastic lightpath rerouting} (\lr{}) problem, for obtaining a provisioning that maximizes the expected number of granted requests in the future. We formulate the \lr{} problem as a 2SP model that shares many similarities with the \srwa{}  model.
	\item Relying on the empirical evidence of its strength, we propose to solve a certain relaxation of the stochastic models to derive high-quality solutions efficiently. For solving the proposed relaxation, we design a decomposition method enhanced by  a new problem-specific family of cuts that significantly reduces the number of optimality cuts needed for its convergence.
	\item Through extensive numerical experiments, we show the notable value of considering stochasticity in provisioning and defragmentation problems in WDM networks, as well as the strength of our new family of cuts. 
\end{itemize}

The rest of the paper is organized as follows. In Section \ref{sec: det-defrag} we briefly review the preliminaries and then introduce two novel stochastic problems to be used in provisioning and defragmentation processes. In Section \ref{sec:problem} we propose 2SP models for these two problems. In Section \ref{sec:benders} we design a decomposition framework for solving the proposed models and enhance it by introducing a new family of cuts. In Section \ref{sec:num-results} we provide a thorough numerical study that evaluates the value of incorporating stochasticity in provisioning and defragmentation, as well as the performance of our  solution method. Finally, in Section \ref{sec:conlusion} we provide concluding remarks and future research directions.  

\section{Background and Motivation}
\label{sec: det-defrag}
Consider a WDM network with a given physical topology and a set of available wavelengths. \rev{The traffic on this network is characterized by the number of connection requests between two nodes of the physical network, and is assumed to be asymmetric (i.e., the links in the network are bidirectional, and the number of requests from node $v_1$ to node $v_2$ can be different than the number of requests from $v_2$ to $v_1$)}. A connection request between two nodes is served by a lightpath, consisting of a path between those nodes together with a certain wavelength. Note that, in the absence of wavelength conversion, shown to be practically ineffective in increasing  GoS \citep{zhang2013proof}, the same wavelength is used throughout the entire path associated with a lightpath. This is known as \emph{wavelength continuity} property. As another important property, two lightpaths cannot share links with the same wavelength, otherwise they incur a \emph{wavelength conflict}. 
\begin{figure}[h]
	\centering
	\small
	\subfigure[Network topology]{%
		\tikzset{every picture/.style={line width=0.75pt}} 

\begin{tikzpicture}[scale=\textwidth/30cm,x=0.75pt,y=0.75pt,yscale=-1,xscale=1]

\draw [line width=0.75]    (134.4,113.55) -- (154.35,90.75) ;
\draw [line width=0.75]    (151.5,137.77) -- (199.95,152.98) ;
\draw [line width=0.75]    (134.4,151.55) -- (163.85,204.75) ;
\draw [line width=0.75]    (172.4,104.05) -- (168.6,201.9) ;
\draw [line width=0.75]    (235.1,162.95) -- (264.55,161.05) ;
\draw [line width=0.75]    (191.4,85.05) -- (396.6,133.5) ;
\draw  [line width=0.75]  (115.4,132.55) .. controls (115.4,122.06) and (123.91,113.55) .. (134.4,113.55) .. controls (144.89,113.55) and (153.4,122.06) .. (153.4,132.55) .. controls (153.4,143.04) and (144.89,151.55) .. (134.4,151.55) .. controls (123.91,151.55) and (115.4,143.04) .. (115.4,132.55) -- cycle ;
\draw  [line width=0.75]  (153.4,85.05) .. controls (153.4,74.56) and (161.91,66.05) .. (172.4,66.05) .. controls (182.89,66.05) and (191.4,74.56) .. (191.4,85.05) .. controls (191.4,95.54) and (182.89,104.05) .. (172.4,104.05) .. controls (161.91,104.05) and (153.4,95.54) .. (153.4,85.05) -- cycle ;
\draw  [line width=0.75]  (149.6,220.9) .. controls (149.6,210.41) and (158.11,201.9) .. (168.6,201.9) .. controls (179.09,201.9) and (187.6,210.41) .. (187.6,220.9) .. controls (187.6,231.39) and (179.09,239.9) .. (168.6,239.9) .. controls (158.11,239.9) and (149.6,231.39) .. (149.6,220.9) -- cycle ;
\draw  [line width=0.75]  (197.1,162.95) .. controls (197.1,152.46) and (205.61,143.95) .. (216.1,143.95) .. controls (226.59,143.95) and (235.1,152.46) .. (235.1,162.95) .. controls (235.1,173.44) and (226.59,181.95) .. (216.1,181.95) .. controls (205.61,181.95) and (197.1,173.44) .. (197.1,162.95) -- cycle ;
\draw  [line width=0.75]  (264.55,161.05) .. controls (264.55,150.56) and (273.06,142.05) .. (283.55,142.05) .. controls (294.04,142.05) and (302.55,150.56) .. (302.55,161.05) .. controls (302.55,171.54) and (294.04,180.05) .. (283.55,180.05) .. controls (273.06,180.05) and (264.55,171.54) .. (264.55,161.05) -- cycle ;
\draw  [line width=0.75]  (271.2,272.2) .. controls (271.2,261.71) and (279.71,253.2) .. (290.2,253.2) .. controls (300.69,253.2) and (309.2,261.71) .. (309.2,272.2) .. controls (309.2,282.69) and (300.69,291.2) .. (290.2,291.2) .. controls (279.71,291.2) and (271.2,282.69) .. (271.2,272.2) -- cycle ;
\draw  [line width=0.75]  (319.65,151.55) .. controls (319.65,141.06) and (328.16,132.55) .. (338.65,132.55) .. controls (349.14,132.55) and (357.65,141.06) .. (357.65,151.55) .. controls (357.65,162.04) and (349.14,170.55) .. (338.65,170.55) .. controls (328.16,170.55) and (319.65,162.04) .. (319.65,151.55) -- cycle ;
\draw  [line width=0.75]  (377.6,152.5) .. controls (377.6,142.01) and (386.11,133.5) .. (396.6,133.5) .. controls (407.09,133.5) and (415.6,142.01) .. (415.6,152.5) .. controls (415.6,162.99) and (407.09,171.5) .. (396.6,171.5) .. controls (386.11,171.5) and (377.6,162.99) .. (377.6,152.5) -- cycle ;
\draw  [line width=0.75]  (407.05,29) .. controls (407.05,18.51) and (415.56,10) .. (426.05,10) .. controls (436.54,10) and (445.05,18.51) .. (445.05,29) .. controls (445.05,39.49) and (436.54,48) .. (426.05,48) .. controls (415.56,48) and (407.05,39.49) .. (407.05,29) -- cycle ;
\draw  [line width=0.75]  (437.45,143) .. controls (437.45,132.51) and (445.96,124) .. (456.45,124) .. controls (466.94,124) and (475.45,132.51) .. (475.45,143) .. controls (475.45,153.49) and (466.94,162) .. (456.45,162) .. controls (445.96,162) and (437.45,153.49) .. (437.45,143) -- cycle ;
\draw  [line width=0.75]  (371.9,276) .. controls (371.9,265.51) and (380.41,257) .. (390.9,257) .. controls (401.39,257) and (409.9,265.51) .. (409.9,276) .. controls (409.9,286.49) and (401.39,295) .. (390.9,295) .. controls (380.41,295) and (371.9,286.49) .. (371.9,276) -- cycle ;
\draw  [line width=0.75]  (522,76.5) .. controls (522,66.01) and (530.51,57.5) .. (541,57.5) .. controls (551.49,57.5) and (560,66.01) .. (560,76.5) .. controls (560,86.99) and (551.49,95.5) .. (541,95.5) .. controls (530.51,95.5) and (522,86.99) .. (522,76.5) -- cycle ;
\draw  [line width=0.75]  (436.5,255.1) .. controls (436.5,244.61) and (445.01,236.1) .. (455.5,236.1) .. controls (465.99,236.1) and (474.5,244.61) .. (474.5,255.1) .. controls (474.5,265.59) and (465.99,274.1) .. (455.5,274.1) .. controls (445.01,274.1) and (436.5,265.59) .. (436.5,255.1) -- cycle ;
\draw  [line width=0.75]  (522,199.05) .. controls (522,188.56) and (530.51,180.05) .. (541,180.05) .. controls (551.49,180.05) and (560,188.56) .. (560,199.05) .. controls (560,209.54) and (551.49,218.05) .. (541,218.05) .. controls (530.51,218.05) and (522,209.54) .. (522,199.05) -- cycle ;
\draw [line width=0.75]    (302.55,161.05) -- (322.5,160.81) ;
\draw [line width=0.75]    (357.18,156.06) -- (378.55,155.11) ;
\draw [line width=0.75]    (415.6,152.5) -- (438.4,147.99) ;
\draw [line width=0.75]    (474.03,135.64) -- (528.17,91.46) ;
\draw [line width=0.75]    (470.23,156.06) -- (523.42,194.54) ;
\draw [line width=0.75]    (228.45,147.28) -- (409.9,39.69) ;
\draw [line width=0.75]    (445.05,29) -- (524.85,66.29) ;
\draw [line width=0.75]    (440.3,40.64) -- (534.35,180.76) ;
\draw [line width=0.75]    (466.43,239.19) -- (537.67,96.21) ;
\draw [line width=0.75]    (453.6,161.76) -- (399.93,258.19) ;
\draw [line width=0.75]    (283.55,180.05) -- (290.2,253.2) ;
\draw [line width=0.75]    (184.28,230.16) -- (273.1,266.74) ;
\draw [line width=0.75]    (372.38,273.39) -- (310.15,274.34) ;
\draw [line width=0.75]    (437.45,250.11) -- (307.78,263.41) ;
\draw [line width=0.75]    (474.5,255.1) -- (535.3,217.58) ;
\draw [line width=0.75]    (144.85,24.25) -- (182.85,24.25) ;
\draw [line width=0.75]    (201.85,14.75) -- (236.85,14.75) ;
\draw [shift={(239.85,14.75)}, rotate = 180] [fill={rgb, 255:red, 0; green, 0; blue, 0 }  ][line width=0.08]  [draw opacity=0] (7.14,-3.43) -- (0,0) -- (7.14,3.43) -- (4.74,0) -- cycle    ;
\draw [line width=0.75]    (204.85,33.75) -- (239.85,33.75) ;
\draw [shift={(201.85,33.75)}, rotate = 0] [fill={rgb, 255:red, 0; green, 0; blue, 0 }  ][line width=0.08]  [draw opacity=0] (7.14,-3.43) -- (0,0) -- (7.14,3.43) -- (4.74,0) -- cycle    ;
\draw  [fill={rgb, 255:red, 0; green, 0; blue, 0 }  ,fill opacity=1 ][line width=0.75]  (192.35,21.08) .. controls (192.35,20.21) and (191.64,19.5) .. (190.77,19.5) .. controls (189.89,19.5) and (189.18,20.21) .. (189.18,21.08) .. controls (189.18,21.96) and (189.89,22.67) .. (190.77,22.67) .. controls (191.64,22.67) and (192.35,21.96) .. (192.35,21.08) -- cycle ;
\draw  [fill={rgb, 255:red, 0; green, 0; blue, 0 }  ,fill opacity=1 ][line width=0.75]  (192.35,27.42) .. controls (192.35,26.54) and (191.64,25.83) .. (190.77,25.83) .. controls (189.89,25.83) and (189.18,26.54) .. (189.18,27.42) .. controls (189.18,28.29) and (189.89,29) .. (190.77,29) .. controls (191.64,29) and (192.35,28.29) .. (192.35,27.42) -- cycle ;

\draw (173.51,85.21) node    {$2$};
\draw (136.16,131.78) node    {$1$};
\draw (170.34,218.84) node    {$3$};
\draw (217.84,161.84) node    {$4$};
\draw (284.34,161.84) node    {$5$};
\draw (290.99,271.41) node    {$6$};
\draw (339.76,152.34) node    {$7$};
\draw (396.76,152.34) node    {$8$};
\draw (426.84,28.84) node    {$9$};
\draw (455.34,142.84) node    {$10$};
\draw (391.06,276.16) node    {$11$};
\draw (540.84,76.34) node    {$12$};
\draw (455.34,255.26) node    {$13$};
\draw (540.84,198.26) node    {$14$};

\end{tikzpicture}
		\label{fig:topology}} 
	\hspace*{.5cm}
	\subfigure[An invalid provisioning]{%
		\input{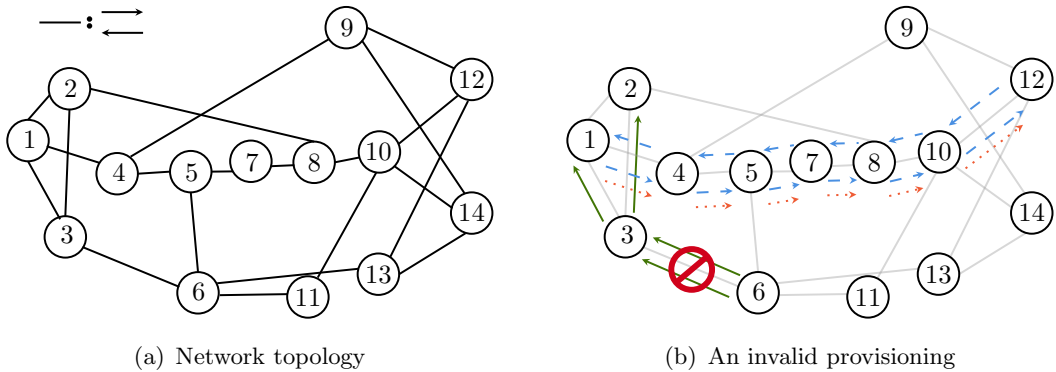}%
		\label{fig:lightpaths}} \\\
	\caption{A WDM network}
	\label{fig:wdm-network}
\end{figure}  

Figure \ref{fig:wdm-network} provides an example of the physical topology of a WDM network, along with a provisioning for two requests from node 1 to 12, one request from node 12 to node 1, one request from node 6 to node 1, and one request from node 6 to node 2. In this figure, every link consists of bidirectional fibers. There are 3 wavelengths available on this network, represented by solid green, dashed blue and dotted red. While dashed and dotted lightpaths respect wavelength continuity with no wavelength conflicts, solid lightpaths share the same wavelength on a fiber link, hence they have a wavelength conflict. A provisioning is called valid if its lightpaths do not violate wavelength continuity and they do not have any wavelength conflicts. 
Thus, the provisioning illustrated in Figure \ref{fig:lightpaths} is not valid. If the request from node 6 to node 1 was granted on the same path with a different wavelength, then the provisioning would have been valid. For a set of incoming connection requests, the \maxRWA{} problem finds a valid provisioning that maximizes the GoS, while \minRWA{} minimizes the number of links used in the provisioning.

\subsection{Dynamic Traffic}
\begin{figure}[htbp]
	\centering
	\small
	\input{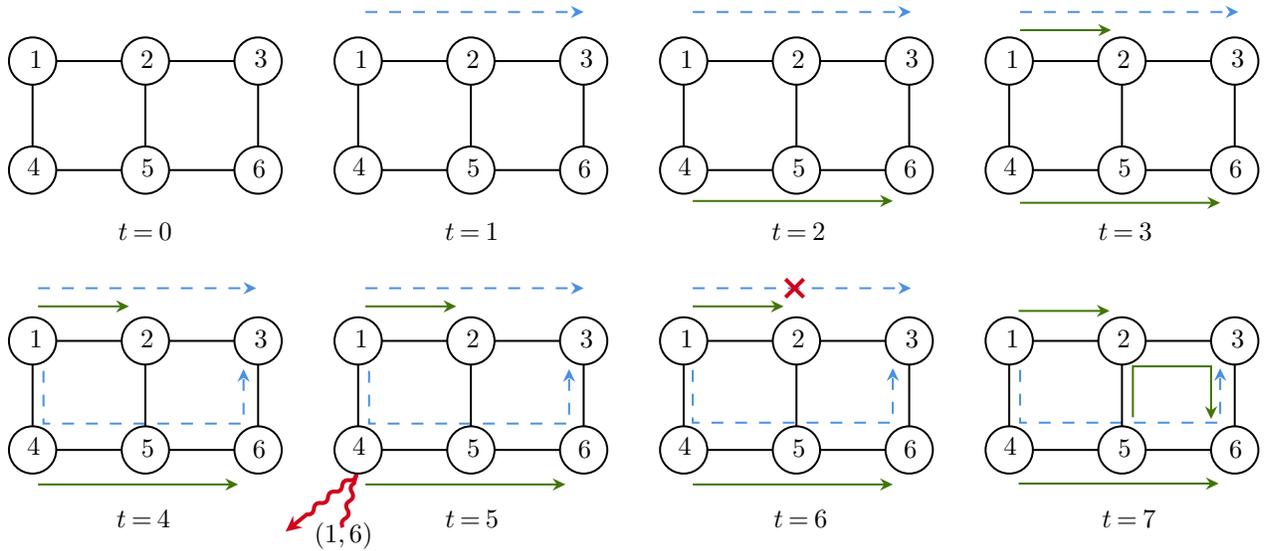}%
	\caption{Dynamic traffic in a WDM network with two wavelengths}
	\label{fig:dynRWA}
\end{figure}

In a WDM network, dynamic traffic involves adds and drops of  connection requests over time, through lightpath assignments and release of wavelinks, respectively. Figure \ref{fig:dynRWA} shows an example of a dynamic traffic in a WDM network with two available wavelengths. Incoming traffic is granted until period $ t = 4 $ (represented by solid green and dashed blue). At $ t=5 $ a request from node 1 to 6 is blocked. At period $ t = 6 $, a connection is leaving the network, i.e., a connection drop occurs, then granting of the connections continues. 

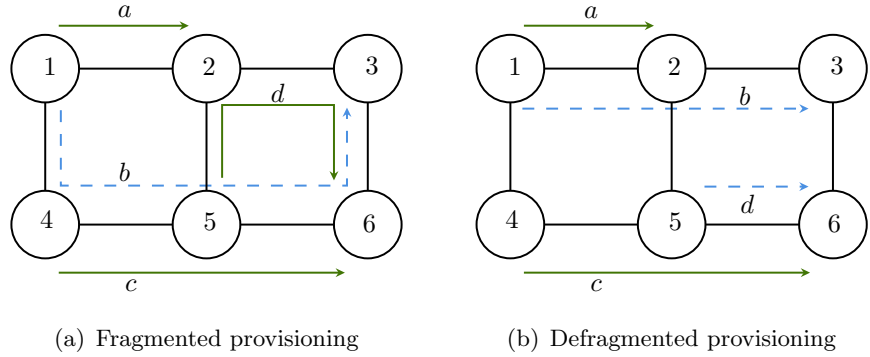
\begin{figure}[t]
	\centering
	\small
	\subfigure[Fragmented provisioning]{%
		\tikzset{every picture/.style={line width=0.75pt}} 

\begin{tikzpicture}[scale=\textwidth/25cm,x=0.75pt,y=0.75pt,yscale=-1,xscale=1]

\draw [color={rgb, 255:red, 65; green, 117; blue, 5 }  ,draw opacity=1 ][line width=0.75]    (216.16,50.07) -- (312.49,50.07) ;
\draw [shift={(315.49,50.07)}, rotate = 180] [fill={rgb, 255:red, 65; green, 117; blue, 5 }  ,fill opacity=1 ][line width=0.08]  [draw opacity=0] (7.14,-3.43) -- (0,0) -- (7.14,3.43) -- (4.74,0) -- cycle    ;
\draw [color={rgb, 255:red, 65; green, 117; blue, 5 }  ,draw opacity=1 ][line width=0.75]    (216.16,238.8) -- (431.69,238.8) ;
\draw [shift={(434.69,238.8)}, rotate = 180] [fill={rgb, 255:red, 65; green, 117; blue, 5 }  ,fill opacity=1 ][line width=0.08]  [draw opacity=0] (7.14,-3.43) -- (0,0) -- (7.14,3.43) -- (4.74,0) -- cycle    ;
\draw [color={rgb, 255:red, 74; green, 144; blue, 226 }  ,draw opacity=1 ][line width=0.75]  [dash pattern={on 4.5pt off 4.5pt}]  (217.75,114.64) -- (217.75,172.25) -- (436.28,172.25) -- (436.28,115.65) ;
\draw [shift={(436.28,112.65)}, rotate = 450] [fill={rgb, 255:red, 74; green, 144; blue, 226 }  ,fill opacity=1 ][line width=0.08]  [draw opacity=0] (7.14,-3.43) -- (0,0) -- (7.14,3.43) -- (4.74,0) -- cycle    ;
\draw [color={rgb, 255:red, 65; green, 117; blue, 5 }  ,draw opacity=1 ][line width=0.75]    (340.92,166.29) -- (340.92,110.66) -- (426.34,110.66) -- (426.34,165.28) ;
\draw [shift={(426.34,168.28)}, rotate = 270] [fill={rgb, 255:red, 65; green, 117; blue, 5 }  ,fill opacity=1 ][line width=0.08]  [draw opacity=0] (7.14,-3.43) -- (0,0) -- (7.14,3.43) -- (4.74,0) -- cycle    ;
\draw   (180,82.85) .. controls (180,68.59) and (191.56,57.03) .. (205.83,57.03) .. controls (220.09,57.03) and (231.65,68.59) .. (231.65,82.85) .. controls (231.65,97.11) and (220.09,108.68) .. (205.83,108.68) .. controls (191.56,108.68) and (180,97.11) .. (180,82.85) -- cycle ;
\draw   (303.17,82.85) .. controls (303.17,68.59) and (314.73,57.03) .. (329,57.03) .. controls (343.26,57.03) and (354.82,68.59) .. (354.82,82.85) .. controls (354.82,97.11) and (343.26,108.68) .. (329,108.68) .. controls (314.73,108.68) and (303.17,97.11) .. (303.17,82.85) -- cycle ;
\draw   (426.34,82.85) .. controls (426.34,68.59) and (437.91,57.03) .. (452.17,57.03) .. controls (466.43,57.03) and (478,68.59) .. (478,82.85) .. controls (478,97.11) and (466.43,108.68) .. (452.17,108.68) .. controls (437.91,108.68) and (426.34,97.11) .. (426.34,82.85) -- cycle ;
\draw   (180,202.05) .. controls (180,187.79) and (191.56,176.22) .. (205.83,176.22) .. controls (220.09,176.22) and (231.65,187.79) .. (231.65,202.05) .. controls (231.65,216.31) and (220.09,227.88) .. (205.83,227.88) .. controls (191.56,227.88) and (180,216.31) .. (180,202.05) -- cycle ;
\draw   (303.17,202.05) .. controls (303.17,187.79) and (314.73,176.22) .. (329,176.22) .. controls (343.26,176.22) and (354.82,187.79) .. (354.82,202.05) .. controls (354.82,216.31) and (343.26,227.88) .. (329,227.88) .. controls (314.73,227.88) and (303.17,216.31) .. (303.17,202.05) -- cycle ;
\draw   (426.34,202.05) .. controls (426.34,187.79) and (437.91,176.22) .. (452.17,176.22) .. controls (466.43,176.22) and (478,187.79) .. (478,202.05) .. controls (478,216.31) and (466.43,227.88) .. (452.17,227.88) .. controls (437.91,227.88) and (426.34,216.31) .. (426.34,202.05) -- cycle ;
\draw    (231.65,82.85) -- (303.17,82.85) ;
\draw    (354.82,82.85) -- (426.34,82.85) ;
\draw    (205.83,108.68) -- (205.83,176.22) ;
\draw    (329,108.68) -- (329,176.22) ;
\draw    (452.17,108.68) -- (452.17,176.22) ;
\draw    (231.65,202.05) -- (303.17,202.05) ;
\draw    (354.82,202.05) -- (426.34,202.05) ;

\draw (453.4,200.86) node    {$6$};
\draw (331.82,200.06) node    {$5$};
\draw (207.06,198.47) node    {$4$};
\draw (455.39,80.47) node    {$3$};
\draw (331.02,80.86) node    {$2$};
\draw (209.44,80.86) node    {$1$};
\draw (266.5,39) node   [align=left] {$\displaystyle a$};
\draw (266.5,162) node   [align=left] {$\displaystyle b$};
\draw (271.37,248.54) node   [align=left] {$\displaystyle c$};
\draw (383.5,101) node   [align=left] {$\displaystyle d$};

\end{tikzpicture}%
		\label{fig:fragmented-example}} 
	\hspace*{.5cm}
	\subfigure[Defragmented provisioning]{%
		\tikzset{every picture/.style={line width=0.75pt}} 

\begin{tikzpicture}[scale=\textwidth/25cm,x=0.75pt,y=0.75pt,yscale=-1,xscale=1]

\draw [color={rgb, 255:red, 65; green, 117; blue, 5 }  ,draw opacity=1 ][line width=0.75]    (214.16,50.53) -- (310.49,50.53) ;
\draw [shift={(313.49,50.53)}, rotate = 180] [fill={rgb, 255:red, 65; green, 117; blue, 5 }  ,fill opacity=1 ][line width=0.08]  [draw opacity=0] (7.14,-3.43) -- (0,0) -- (7.14,3.43) -- (4.74,0) -- cycle    ;
\draw [color={rgb, 255:red, 65; green, 117; blue, 5 }  ,draw opacity=1 ][line width=0.75]    (214.16,239.27) -- (429.69,239.27) ;
\draw [shift={(432.69,239.27)}, rotate = 180] [fill={rgb, 255:red, 65; green, 117; blue, 5 }  ,fill opacity=1 ][line width=0.08]  [draw opacity=0] (7.14,-3.43) -- (0,0) -- (7.14,3.43) -- (4.74,0) -- cycle    ;
\draw [color={rgb, 255:red, 74; green, 144; blue, 226 }  ,draw opacity=1 ][line width=0.75]  [dash pattern={on 4.5pt off 4.5pt}]  (213.01,113.91) -- (428.54,113.91) ;
\draw [shift={(431.54,113.91)}, rotate = 180] [fill={rgb, 255:red, 74; green, 144; blue, 226 }  ,fill opacity=1 ][line width=0.08]  [draw opacity=0] (7.14,-3.43) -- (0,0) -- (7.14,3.43) -- (4.74,0) -- cycle    ;
\draw [color={rgb, 255:red, 74; green, 144; blue, 226 }  ,draw opacity=1 ][line width=0.75]  [dash pattern={on 4.5pt off 4.5pt}]  (352.07,173.51) -- (428.54,173.51) ;
\draw [shift={(431.54,173.51)}, rotate = 180] [fill={rgb, 255:red, 74; green, 144; blue, 226 }  ,fill opacity=1 ][line width=0.08]  [draw opacity=0] (7.14,-3.43) -- (0,0) -- (7.14,3.43) -- (4.74,0) -- cycle    ;
\draw   (178,83.31) .. controls (178,69.05) and (189.57,57.49) .. (203.83,57.49) .. controls (218.09,57.49) and (229.66,69.05) .. (229.66,83.31) .. controls (229.66,97.58) and (218.09,109.14) .. (203.83,109.14) .. controls (189.57,109.14) and (178,97.58) .. (178,83.31) -- cycle ;
\draw   (301.18,83.31) .. controls (301.18,69.05) and (312.74,57.49) .. (327,57.49) .. controls (341.27,57.49) and (352.83,69.05) .. (352.83,83.31) .. controls (352.83,97.58) and (341.27,109.14) .. (327,109.14) .. controls (312.74,109.14) and (301.18,97.58) .. (301.18,83.31) -- cycle ;
\draw   (424.35,83.31) .. controls (424.35,69.05) and (435.91,57.49) .. (450.17,57.49) .. controls (464.44,57.49) and (476,69.05) .. (476,83.31) .. controls (476,97.58) and (464.44,109.14) .. (450.17,109.14) .. controls (435.91,109.14) and (424.35,97.58) .. (424.35,83.31) -- cycle ;
\draw   (178,202.51) .. controls (178,188.25) and (189.57,176.69) .. (203.83,176.69) .. controls (218.09,176.69) and (229.66,188.25) .. (229.66,202.51) .. controls (229.66,216.78) and (218.09,228.34) .. (203.83,228.34) .. controls (189.57,228.34) and (178,216.78) .. (178,202.51) -- cycle ;
\draw   (301.18,202.51) .. controls (301.18,188.25) and (312.74,176.69) .. (327,176.69) .. controls (341.27,176.69) and (352.83,188.25) .. (352.83,202.51) .. controls (352.83,216.78) and (341.27,228.34) .. (327,228.34) .. controls (312.74,228.34) and (301.18,216.78) .. (301.18,202.51) -- cycle ;
\draw   (424.35,202.51) .. controls (424.35,188.25) and (435.91,176.69) .. (450.17,176.69) .. controls (464.44,176.69) and (476,188.25) .. (476,202.51) .. controls (476,216.78) and (464.44,228.34) .. (450.17,228.34) .. controls (435.91,228.34) and (424.35,216.78) .. (424.35,202.51) -- cycle ;
\draw    (229.66,83.31) -- (301.18,83.31) ;
\draw    (352.83,83.31) -- (424.35,83.31) ;
\draw    (203.83,109.14) -- (203.83,176.69) ;
\draw    (327,109.14) -- (327,176.69) ;
\draw    (450.17,109.14) -- (450.17,176.69) ;
\draw    (229.66,202.51) -- (301.18,202.51) ;
\draw    (352.83,202.51) -- (424.35,202.51) ;

\draw (265.4,40.46) node   [align=left] {$\displaystyle a$};
\draw (382.86,104.01) node   [align=left] {$\displaystyle b$};
\draw (269.38,249) node   [align=left] {$\displaystyle c$};
\draw (384.85,186.41) node   [align=left] {$\displaystyle d$};
\draw (451.4,201.32) node    {$6$};
\draw (329.82,200.53) node    {$5$};
\draw (205.06,198.94) node    {$4$};
\draw (453.39,80.93) node    {$3$};
\draw (329.03,81.33) node    {$2$};
\draw (207.44,81.33) node    {$1$};

\end{tikzpicture}%
		\label{fig:defragmented-example}} \\\
	\caption{Defragmentation example}
	\label{fig:defrag-example}
\end{figure}

The consecutive adds and drops over time can result in a fragmented network that is probable to higher rate of request rejections, while there exists capacity over the network. The process of migrating a fragmented network to a provisioning with recovered bandwidth is called defragmentation. Figure \ref{fig:defrag-example} presents a fragmented provisioning in the WDM network of Figure \ref{fig:dynRWA} at $ t=7 $, as well as a defragmented one that grants the same requests while using a fewer number of wavelinks. In a fragmented network, lightpaths with longer lengths are employed in the provisioning, leading to higher wavelink usage.  For example, in Figures \ref{fig:fragmented-example} and \ref{fig:defragmented-example} the average lightpath length is 2.5 and 1.5 links, respectively. The provisioning in Figure \ref{fig:defragmented-example} has recovered the stranded bandwidth by using shorter lightpaths. For instance, the defragmented network can grant four more requests from node 2 to 5, while the fragmented one has the capacity of only three more requests between the same nodes.

\subsection{Motivation}
A deployed and activated WDM network starts with full capacity. Over time, it grants the incoming connection requests by establishing lightpaths and frees the resources once  a connection leaves the network. At certain points,  we want to reroute the current existing lightpaths such that the network has a higher chance of granting the uncertain future demand. From this description, stochastic traffic impacts the performance of a network from two aspects: provisioning and defragmentation.

For the provisioning problem, disregarding the future fosters greedy assignment of lightpaths. It aims for locally maximum number of granted requests, even if it leads to selection of long lightpaths that take up considerable amount of resources just to accept one more connection at the current stage. By considering the future traffic one can make more informed lightpath assignment decisions, both in terms of acceptance or rejection of requests, and the characteristics of lightpaths that are assigned to the granted ones. Consequently, the goal is not to find locally optimal solutions at every stage, rather it is to optimize the overall performance in the longer term.

As for the defragmentation problem, there are two main steps in the process: ($i$) design of a defragmented provisioning target for the current set of connection requests, and ($ii$) conversion of the current provisioning to the target, namely performing a \emph{network migration}. Common provisioning target in the network defragmentation problem is to minimize the usage of resources (e.g., minimization of the number of used links), with the assumption that it leads to a more efficient answer to traffic requests, leaving a larger number of wavelinks to serve the future traffic.Note that, over a multi-period planning horizon, merely minimizing the number of used resources as a local decision does not guarantee the best performance in the future, and despite the momentary decrease in the wavelink usage, the network can quickly become fragmented once again after defragmentation.
In fact, to the best of our knowledge, the impact of the defragmentation outcome on future GoS has not been studied before. 
We argue that a deterministic approach does not have the means to recognize and free up critical wavelinks that can contribute to granting more future requests, thus  such a  provisioning, although defragmented, might not be able to improve the GoS overall.

Following the discussions above, we suggest to look one step ahead while making the lightpath assignment/rerouting decisions, by solving two stochastic programming problems that correspond to stochastic versions of the \maxRWA{} and lightpath rerouting problems, i.e., \srwa{} and \lr{}, respectively. 
Throughout the paper, we make the following set of assumptions:
\begin{enumerate}[label=(A\arabic*), leftmargin=*]
	\item Lightpath assignment decisions are made periodically, specifically all the requests arrived in a time interval are granted/rejected at the end of the interval.
	\item The target of the defragmentation process is reachable, i.e., we can reroute the lightpaths without service interruption.\label{assumprion2}
\end{enumerate}
The first assumption is justified by the common practice of service providers, where the granting decisions are made for a batch of requests in order to maximize GoS and  incremental planning has to be performed in given periods to determine the capacity upgrades needed for addressing the growth of the traffic. It is indeed a common practice in the context of the dynamic traffic literature \citep{wu2012forward,chen2015spectrum,jaumard2017efficient}. The second assumption is nonrestrictive, since by leasing one unit of wavelength \citep{leung2005capacity}, any provisioning can be migrated to a new one without service disruption. 

In the following sections, we formally define \srwa{} and \lr{} and formulate them as 2SP models.

\section{Problem Descriptions and Formulations}
\label{sec:problem}

The stochastic \maxRWA{} and lightpath rerouting problems, namely \srwa{} and \lr{}, both assign lightpaths to a known set of connection requests, such that the expected GoS is maximized. Therefore, their mathematical models are quite similar, with the exception of the definition of node pair set for which they are building a provisioning. \srwa{} attempts to choose from the set of new connection requests, while \lr{} has to assign a (potentially new) lightpath to every connection already on the network. In other words, defragmentation (\lr{}) is performed at the beginning of a time period before the arrival of a new batch of connection requests, while \srwa{} is solved at the end of a period, when the new traffic has arrived, and future traffic is taken into consideration in both. In the following, we first present the 2SP formulation of \srwa{} in Section \ref{sec:SmaxRWA}, then in Section \ref{sec:SmaxLR} we make the necessary modifications for \lr{}.

\subsection{The Stochastic \maxRWA{} Problem}\label{sec:SmaxRWA}

Consider a WDM network. Denote by $\Graph = (\nodeSet,\edgeSet) $ the multigraph that represents the physical topology of this WDM network, with $ \node \in \nodeSet $  and $ \edge \in \edgeSet $ a node and fiber link of the physical layer, respectively. $ \wavelength\in \waveSet $ is an available wavelength on the network. To grant a connection request between a node pair $ (s,d)\in\nodeSet\times\nodeSet $, a
lightpath $ \lightpath = (\LPpath,\wavelength) $ with path  $ \LPpath  $ between $ s $ and $ d $ over a wavelength $ \wavelength $ is set up. In the following, the uncertainty in the future traffic is characterized by $ \Dnew $, denoting the number of requested lightpaths between node pair $(s,d)$, where $ \xirand $ is the underlying random vector with known probability distribution $ \P $ and support $ \Xi $. \rev{The vector $\xirand$ contains any exogenous element impacting the number of requests between pairs of nodes. Then, random variables $\Dnew$ are functions of this input that map the underlying uncertainty to the traffic demand. In the simplest form, $\xirand$ can be the vector consisting of the random demands.}
Suppose that a new batch of connection requests has arrived on a potentially non-empty network. Available wavelinks are denoted by $ \edgeSetw $, the set of links without any granted requests on wavelength $ \wavelength$. 
    \begin{table}[htbp]
    \centering
    \caption{\rev{Notation used in the \srwa{} model}}
    {\begin{tabular}{lp{14.5cm}}
    \toprule
    \multicolumn{2}{l}{\it Sets and parameters:} \\
        $\Graph$& A multigraph representing the physical topology of the WDM network\\
        $\nodeSet$& Set of nodes in $\Graph$, indexed by $\node$\\
        $\edgeSet$& Set of links in $\Graph$, indexed by  $\edge$\\
        $\waveSet$ & Set of available wavelengths, indexed by  $\wavelength$ \\
        $\LPSet$ & Set of available lightpaths, indexed by $\lightpath$\\
        $ \edgeSetw $ & The set of free links on wavelength $\wavelength$ \\
        $\nodefunc^{+ / -}(\node)$& Set of links leaving/entering node $\node$\\
        $\xirand$& Vector of random variables, with probability distribution $\P$ and support $\xiset$ \\
        $\SD$& Set of node pairs with at least one request
        in $\Graph$, indexed by $(s,d)$
        \\
        $\SDcurr{n}$ & Set of node pairs with at least one new request\\
        $\SDcurr{c}$ &Set of node pairs with at least one existing request\\
        $\SDsce$& Set of node pairs with at least one  request in the future\\
        $\Dcurr$& Number of current connections requests for node pair $(s,d)$ \\
        $\Dnew$& Number of future connections requests for node pair $(s,d)$ \\
        \midrule
        \multicolumn{2}{l}{\it Decision variables:} \\
        $\link$& 1 if link $\edge$ on wavelength $\wavelength$ is serving node pair $(s,d)$, 0 otherwise \\
        $\linksce$& 1 if link $\edge$ on wavelength $\wavelength$ is serving  node pair $(s,d)$ in the future, 0 otherwise  \\
        \bottomrule
    \end{tabular}
    }
    \label{tab:parameters}
\end{table}

The immediate (first-stage) decisions \srwa{} has to make are the lightpath assignments for the new incoming traffic. Therefore the first-stage problem of \srwa{} is a \maxRWA{} problem with a slightly different objective function, given a set of available wavelinks. Let  $\SDcurr{n}\subseteq\{(s,d) \in \nodeSet\times \nodeSet\}$ be the set of node pairs with at least one new connection request in the most recent batch of arrivals, ($ \textsc{n} $ stands for ``new'').  A node pair $ (s, d) $ can have multiple connection requests, denoted by \rev{$\Dcurr\in\N$, with $\N$ the set of natural numbers}.   $\link$ is a binary decision variable equal to 1 if link $\edge$ is used in a lightpath with wavelength $\wavelength$ for a connection request between node pair $(s,d)$. \rev{$x$ is the vector of these decision variables.}
\rev{Table \ref{tab:parameters} summarizes all the notation used in the model.}

Given $ \SDcurr{n} $, \srwa{} is formulated as follows:
\BSE
\label{eq: SmaxRWA}
\begin{alignat}{5}
\max \ \ & \sum_{(s,d) \in \SDcurr{n}}\sum_{\wavelength \in \waveSet}\sum_{\stackrel{\edge\in}{\nodefunc^+(s)\cap \edgeSetw}}\link+\E_{\xirand}\big[\Q(x,\xirand,\SDcurr{n})\big] \label{eq: OF}\\
\text{s.t.}\ \ & x\in \Feas{\SDcurr{n}},
\end{alignat}
\ESE
where $ \Q(.) $ is  the value function that measures GoS in the second stage and $ \Feas{\SDcurr{n}} $ is the feasible set of the \srwa{} problem. Objective function \eqref{eq: OF} maximizes the expected number of granted requests, considering both the new connection requests and future traffic. $ \Q(.) $ and $ \Feas{.} $ are also used in the modeling  of \lr{}, so in the following we define them for a  generic set of node pairs to avoid repetition. For a given $ \SD $, $ \Feas{\SD} $ is defined by the set of following constraints:
\BSE
\label{eq: Feas}
\begin{alignat}{5}
	& \sum_{(s,d) \in \SD} \link \leq 1 &\quad & \wavelength\in \waveSet, \edge \in \edgeSetw \label{eq: conflict_link}\\
	&\sum_{\stackrel{\edge \in}{\nodefunc^+(\node)\cap \edgeSetw}} \link = \sum_{\stackrel{\edge \in}{ \nodefunc^-(\node)\cap \edgeSetw}} \link && \wavelength \in \waveSet, (s,d) \in \SD,\node\in \nodeSet\setminus\{s,d\} \label{eq: flow_conservation}\\
	&\sum_{\wavelength \in \waveSet} \sum_{\stackrel{\edge \in}{\nodefunc^-(s)\cap\edgeSetw}}\link =\sum_{\wavelength \in \waveSet}  \sum_{\stackrel{\edge \in}{\nodefunc^+(d)\cap\edgeSetw}} \link = 0 && (s,d) \in \SD \label{eq: loop_elimination}\\
	&\sum_{\wavelength \in \waveSet} \sum_{\stackrel{\edge \in}{\nodefunc^+(s)\cap \edgeSetw}} \link \leq \Dcurr && (s,d) \in \SD \label{eq: demand_UB} \\
	& \link\in\{0,1\} && (s,d) \in \SD, \wavelength\in\waveSet,\edge\in\edgeSetw,\label{eq:links_bounds}
\end{alignat}
\ESE
where $\nodefunc^{+(-)}(\node)$ is the set of links leaving (entering) node $\node\in \nodeSet$. 
Constraints \eqref{eq: conflict_link} ensure that every link $\edge$ can be used for at most one path in any wavelength for the new traffic.
Constraints \eqref{eq: flow_conservation} establish the flow conservation in the chosen lightpaths. They  make sure that for a given node pair $ (s,d) $, the same number of wavelinks over the same wavelength is leaving and entering the nodes of the network for granting a request of $ (s,d) $. Constraints \eqref{eq: loop_elimination} eliminate the cycles by preventing the selection of incoming/outgoing link to/from source/destination of a request, and together with constraints \eqref{eq: flow_conservation} set up a proper lightpath for a granted request.
Constraints \eqref{eq: demand_UB} bound the number of lightpath assignments to a node pair by its demand.
Constraints \eqref{eq:links_bounds} are the integrality and bound constraints.
 
With the first-stage provisioning fixed, recourse decisions are made by solving another \maxRWA{} problem for the future traffic $ \Dnew $ on the available resources. $\linksce$ is a second-stage decision variable equal to 1 if for a future connection request between $ (s,d)\in \SDsce=\{(s,d) \in \nodeSet\times \nodeSet: \Dnew > 0\}$, a wavelink with link $ \edge $ and wavelength $ \wavelength $ is used.
$\numgrantedsce$ is an integer (auxiliary) recourse decision variable measuring the number of granted requests between node pair $(s,d)$. For a first-stage solution $ \hat{x} $ and a future traffic scenario $ \xirand $, the recourse problem is defined as:\\
$\Q(x,\xirand,\SD) =$
\BSE
\begin{alignat}{5} 
\max \ \ & \sum_{(s,d) \in \SDsce} \numgrantedsce \label{eq: second_stage_obj}\\
\text{s.t.}\ \ &\sum_{(s,d) \in \SDsce} \linksce  \leq 1- \sum_{(s,d) \in \SD} \link &\ \ & \wavelength\in \waveSet, \edge \in \edgeSetw \label{eq: conflict_link_demands_scenarios}\\
&\sum_{\stackrel{\edge \in}{\nodefunc^+(\node)\cap \edgeSetw}} \linksce = \sum_{\stackrel{\edge \in}{\nodefunc^-(\node)\cap \edgeSetw}} \linksce && \wavelength \in \waveSet, (s,d) \in \SDsce, \node\in \nodeSet\setminus\{s,d\}\label{eq: flow_conservation_scenarios} \\
& \sum_{\wavelength \in \waveSet} \sum_{\stackrel{\edge \in}{\nodefunc^-(s)\cap \edgeSetw}}\linksce  =\sum_{\wavelength \in \waveSet} \sum_{\stackrel{\edge \in}{\nodefunc^+(d)\cap \edgeSetw}}  \linksce = 0 && (s,d) \in \SDsce\label{eq: loop_elimination_scenarios}\\
&\numgrantedsce \leq \Dnew && (s,d) \in \SDsce \label{eq: demand_UB_scenarios} \\
&\numgrantedsce = \sum_{\wavelength \in \waveSet}\sum_{\stackrel{\edge \in}{\nodefunc^+(s)\cap \edgeSetw}} \linksce && (s,d) \in \SDsce\label{eq: granted_scenarios}\\
& \linksce \in\{0,1\} && (s,d) \in \SDsce,\wavelength\in\waveSet, \edge\in\edgeSetw\label{eq: linksce_bounds}\\
&\numgrantedsce\in\Z_+ && (s,d) \in \SDsce.\label{eq: granted_scenarios_bounds}
\end{alignat}
\ESE
Objective function \eqref{eq: second_stage_obj} maximizes the number of granted requests. 
If a wavelink is used in the first-stage provisioning, Constraints	\eqref{eq: conflict_link_demands_scenarios} prevent it to be used in future lightpath assignments.
Constraints \eqref{eq: flow_conservation_scenarios}-\eqref{eq: demand_UB_scenarios}, are the counterparts of Constraints \eqref{eq: flow_conservation}-\eqref{eq: demand_UB} for the second-stage provisioning, ensuring that a proper lightpath is selected for each granted request.
Constraints \eqref{eq: granted_scenarios} calculate the number of granted connections for node pairs, for each realization of future traffic. Constraints \eqref{eq: linksce_bounds}-\eqref{eq: granted_scenarios_bounds} are the integrality and bound constraints in the second stage.

\subsection{Stochastic Lightpath Rerouting Problem}\label{sec:SmaxLR}
In the defragmentation problem, to decide on lightpath reroutings, thanks to assumption \ref{assumprion2} it is irrelevant which wavelinks have been currently in use for the current provisioning. Therefore, the first-stage decisions simply choose the wavelinks in a new provisioning for the current set of granted requests on the network. In other words, in the first-stage problem we have $ \edgeSetw = \edgeSet, \wavelength\in\waveSet $. 
After the defragmentation procedure is complete and the uncertainty is revealed, a \maxRWA{} problem is solved to grant the new connection requests. Considering this,  we propose a 2SP model for the \lr{} problem, where in the first stage we make the rerouting decisions and in the second stage, with the occupied wavelinks now as given, we solve a \maxRWA{} problem for the future traffic. Note that, the definition of the second-stage problem is the same as the recourse problem of \srwa{}.

Let $\SDcurr{c}$ be the set of node pairs with an existing connection on the network ($ \textsc{c} $ stands for ``current''). With the definitions and notations given in Section \ref{sec:SmaxRWA}, the 2SP formulation for \lr{} problem is:
\BSE
\label{eq: SmaxLR}
\begin{alignat}{5}
\max \ \ & \E_{\xirand}\big[\Q(x,\xirand,\SDcurr{c})\big] \label{eq: OF-SmaxLR}\\
\text{s.t.}\ \ & x\in \Feas{\SDcurr{c}}\label{eq: common-SmaxLR}\\
&\sum_{\wavelength \in \waveSet} \sum_{\edge \in\nodefunc^+(s)} \link = \Dcurr &\quad& (s,d) \in \SDcurr{c}. \label{eq: demand}
\end{alignat}
\ESE

As was mentioned earlier, \srwa{} and \lr{} formulations share many commonalities, such as constraints \eqref{eq: common-SmaxLR} and their recourse problem. Differently in \lr{}, every connection already on the network needs to be granted in the new defragmented provisioning, which is modeled by constraints \eqref{eq: demand}, and the objective function \eqref{eq: OF-SmaxLR} measures the expected GoS for the future traffic. Note that constraints \eqref{eq: demand_UB_scenarios} are redundant for the \lr{} model due to constraints  \eqref{eq: demand}.

In the next section, we discuss how the 2SP formulation of \srwa{} and \lr{} problems can be solved efficiently.

\section{Solution Methodology}
\label{sec:benders}
In this section, we use \emph{sample average approximation} (SAA),  to obtain a deterministic equivalent (DE) form of models  \eqref{eq: SmaxRWA} and \eqref{eq: SmaxLR} for a given sample and discuss how to solve it by means of decomposition. Then we introduce a new family of cuts that enhance the performance of the decomposition algorithm. 
\subsection{Sample Average Approximation}\label{sec:saa}
The expectation term in the objective functions of \srwa{} and \lr{} models includes integration over a (high dimensional) random vector $ \xirand $. To overcome the difficulty of maximizing such a function, we use SAA, a technique common in solving such problems (see for e.g., \cite{shapiro2014lectures}). Given a sample of future traffic scenarios, SAA replaces the expectation in the objective function with the sample average. 
Assume that $ \xisetSAA\subseteq \xiset $ is a set of independent and identically distributed (i.i.d.) observations drawn from $ \P $ using Monte Carlo sampling. For \srwa{}, SAA gives the following approximation of  \eqref{eq: SmaxRWA}:
\BSE
\label{eq: saa-SmaxRWA}
\begin{alignat}{5}
\max \ \ & \sum_{(s,d) \in \SDcurr{n}}\sum_{\wavelength \in \waveSet}\sum_{\stackrel{\edge\in}{\nodefunc^+(s)\cap \edgeSetw}}\link + \frac{1}{|\xisetSAA|} \sum_{\xival \in \xisetSAA} \eta_\xival \label{eq: OF_SAA}\\
\text{s.t.}\ \ & x\in \Feas{\SDcurr{n}} \\
& \eta_\xival \leq \Q(x,\xival,\SDcurr{n}) &\qquad& \xival \in \xisetSAA\\
& \eta_\xival \geq 0 && \xival \in \xisetSAA, \label{eq: eta_bound_SAA}
\end{alignat}
\ESE
where $ \eta_\xival $ is a continuous decision variable representing the recourse value function under scenario $ \xival $. As the sample size increases, the solution to  \eqref{eq: saa-SmaxRWA} converges to that of the original problem (under certain mild conditions, as discussed in \cite{shapiro1996simulation}).  By creating copies of the second-stage decision variables and constraints for each scenario in the sample, we can reformulate the SAA problem \eqref{eq: saa-SmaxRWA} as a mixed-integer program, called  the \emph{extensive form} (given in 
\revacc{{\App}, Section \ref{appendix:EXT-form}}). 
The extensive form  often does not scale well with the number of scenarios, thus we  try to solve it by a decomposition method.

The presence of integer decision variables in the second stage makes  \eqref{eq: saa-SmaxRWA} quite challenging. The logic-based 
Benders decomposition \citep{hooker2003logic}, in particular the integer L-shaped method as its special case \citep{laporte1993integer, angulo2016improving}, is a common solution method for solving the DE obtained from applying SAA on stochastic integer programs. However, such methods are often not efficient. 
As was mentioned in Section \ref{sec:SmaxRWA}, our recourse function is itself a \maxRWA{} problem on available wavelinks. In the literature of the RWA problem with asymmetrical traffic, it has been observed that the linear programming (LP) relaxation of its commonly used IP formulation tends to give a high-quality upper bound \citep{jaumard2007comparison, christodoulopoulos2010offline}. Based on this observation (also confirmed by our extensive experiments in Section \ref{sec:num-results}), we consider a relaxation of our DE where second-stage decisions are relaxed, i.e., assumed to be continuous.
We call this relaxation \IPLP{}, in contrast to \IPIP{} as the original exact IP model.  An important advantage of \IPLP{} is that, despite being a relaxation, it provides feasible first-stage decisions which can be directly implemented in practice. It can also provide a lower bound by an out-of-sample evaluation of its first-stage solution. 

\subsection{Decomposition Algorithm}
\IPLP{} is a mixed-binary linear program with binary first-stage and continuous recourse decision variables, for which Bender decomposition \citep{benders2005partitioning} is an established solution method. It decomposes the problem into a master problem and a subproblem that decide on the first- and second-stage variables, respectively. For 2SPs, the latter further decomposes by scenario, thus leads to computational improvement. Recourse value function $ \Q(.) $ is iteratively approximated by solving the subproblems that return feasibility and optimality cuts, until no more cuts are found or a given optimality gap is achieved. Our second-stage feasibility space is  nonempty for feasible first-stage decisions, since granting none of the future requests is always a feasible solution. Therefore, regardless of the $ \hat{x} $ values, as long as they are feasible, second stage remains feasible, i.e., we have relatively complete recourse, thus optimality cuts are solely needed. The Benders master problem (the so-called multi-cut version) for \IPLP{} is:
\BSE
\begin{alignat}{5}
\master(\xisetSAA)\ : \	\max \ \ & \sum_{(s,d) \in \SD}\sum_{\wavelength \in \waveSet}\sum_{\stackrel{\edge\in}{\nodefunc^+(s)\cap \edgeSetw}}\link+\frac{1}{|\xisetSAA|} \sum_{\xival \in \xisetSAA} \eta_\xival \label{eq: OF_SAA_master}\\
	\text{s.t.}\ \ & \eqref{eq: conflict_link} - \eqref{eq:links_bounds}, \eqref{eq: eta_bound_SAA} \\
	& (\eta_\xival,x) \in \cuts &\qquad& \xival \in \xisetSAA, \label{eq:optimality-cuts}
\end{alignat}
\ESE
where $ \cuts $ is the set described by optimality cuts, which we refer to as \emph{$ x- $cuts}, for scenario $ \xival $.
For a given scenario $ \xival\in\xisetSAA $ and a first-stage solution $\hat{x}$, the Benders subproblem is:
\\*[0.2cm]
\noindent $\sub(\hat{x},\xival):$
\BSE
\begin{alignat}{3}
\max \ & \sum_{(s,d) \in \SDsceSAA} \numgrantedsceSAA\label{eq: OF_bd_sp_beta}\\
\text{s.t.} \ &\sum_{(s,d) \in \SDsceSAA} \linksceSAA\leq 1- \sum_{(s,d) \in \SDcurr{n}} \linkhat  && \wavelength\in \waveSet, \edge \in \edgeSetw \label{eq: conflict_link_demands_scenarios_bd_sp_x_beta}\\
&\sum_{\stackrel{\edge \in}{\nodefunc^{+}(\node)\cap\edgeSetw}} \linksceSAA  = \sum_{\stackrel{\edge \in}{\nodefunc^{-}(\node)\cap\edgeSetw}} \linksceSAA \qquad  && \wavelength \in \waveSet, (s,d) \in \SDsceSAA,v\in V\setminus\{s,d\} \label{eq: flow_conservation_scenarios_bd_sp_beta} \\
&\sum_{\wavelength \in \waveSet}\sum_{\stackrel{\edge \in}{\nodefunc^{-}(s)\cap\edgeSetw}} \linksceSAA  =\sum_{\wavelength \in \waveSet} \sum_{\stackrel{\edge \in}{\nodefunc^{+}(d)\cap\edgeSetw}} \linksceSAA = 0 \ \ && (s,d) \in \SDsceSAA\label{eq: loop_elimination_scenarios_bd_sp_beta}\\
&\sum_{\wavelength \in \waveSet}\sum_{\stackrel{\edge \in}{ \nodefunc^{+}(s)\cap\edgeSetw}} \linksceSAA \leq \DnewSAA &&(s,d) \in \SDsceSAA \label{eq: demand_UB_scenarios_bd_sp_beta} \\
&\numgrantedsceSAA = \sum_{\wavelength \in \waveSet}\sum_{\stackrel{\edge \in}{ \nodefunc^+(s)\cap\edgeSetw}} \linksceSAA && (s,d) \in \SDsceSAA\label{eq: granted_scenarios_sub}\\
&\linksceSAA \geq 0  && (s,d) \in \SDsceSAA, \wavelength \in \waveSet, \edge \in \edgeSetw. \label{eq:granted_sub_domain}
\end{alignat}
\ESE
Let $(\hat{x},\hat{\eta})$ be an optimal solution of the master problem, and $\hat{z}$ contain the optimal $z$ variable values of the subproblem $\sub(\hat{x},\xival)$ for scenario $ \xival \in \xisetSAA$. If $ \hat{\eta}_\xival >   \sum\limits_{(s,d) \in \SDsceSAA} \numgrantedsceSAAhat$, an optimality cut ($ x $-cut) in the form of
$$ \eta_\xival \leq   \sum_{\wavelength \in \waveSet}\sum_{\edge \in \edgeSetw} \dualhat^{\eqref{eq: conflict_link_demands_scenarios_bd_sp_x_beta}}_{\edge,\wavelength} \Big(1-\sum_{(s,d) \in \SD} \link\Big) + \sum_{(s,d) \in \SDsceSAA} \dualhat^{\eqref{eq: demand_UB_scenarios_bd_sp_beta}}_{(s,d)}\DnewSAA$$
is added to the description of $ \cuts $, 
where $ \dualhat^{\eqref{eq: conflict_link_demands_scenarios_bd_sp_x_beta}} $ and $ \dualhat^{\eqref{eq: demand_UB_scenarios_bd_sp_beta}} $ are the optimal dual solutions associated with constraints \eqref{eq: conflict_link_demands_scenarios_bd_sp_x_beta} and \eqref{eq: demand_UB_scenarios_bd_sp_beta}, respectively.

SAA problem and optimality cuts for \lr{} are very similar to the ones for \srwa{}.  The only modifications (as in their 2SP formulations) are the lack of first-stage costs, and respectively using node pair set $ \SDcurr{c} $ and constraints \eqref{eq: demand} instead of $ \SDcurr{n} $ and \eqref{eq: demand_UB}.

In Section \ref{sec:num-results}, we evaluate the performance of the described Benders decomposition algorithm and show that, although the relaxation returns high-quality solutions, solving it with the above mentioned standard framework (even with some well-known enhancements such as implementing it in a branch-and-cut framework) is still computationally difficult. In the following we attempt to overcome this by introducing a new family of problem-specific cuts that improves the scalability of our decomposition algorithm.

\subsection{A New Family of Cuts}
The discussions in this section apply to \srwa{}. For \lr{}, the same results can be obtained via modifications discussed in the previous subsection.  
We define for each $\edge \in \edgeSet$ a new decision variable $ \beta_\edge $, measuring the total number of times link $ \edge $ is used for serving traffic requests in the first stage, and add the constraints
\BSE
\begin{alignat}{5}
& \sum_{\wavelength \in \waveSet} \sum_{(s,d) \in \SDcurr{n}} \link = \beta_\edge &\qquad& \edge \in \edgeSet\\
& (\eta_\xival,\beta) \in \betacuts && \xival\in\xisetSAA \label{eq:betacuts_set}\\
&\beta_\edge \geq 0 && \edge \in \edgeSet,
\end{alignat}
\ESE
to the master problem, resulting in a problem denoted by $ \overline{\master}(\xisetSAA) $. In constraints \eqref{eq:betacuts_set}, $ \betacuts $ is the set described by the inequalities called \emph{$ \beta- $cuts} of the following form:
$$ \eta_\xival \leq  \sum_{\edge \in \edgeSet} \dualhat^{\eqref{eq: conflict_link_demands_scenarios_bd_sp_beta}}_{\edge} (|\edgeSetw|-\beta_\edge) + \sum_{(s,d) \in \SDsceSAA} \dualhat^{\eqref{eq: demand_UB_scenarios_bd_sp_beta_2}}_{(s,d)}\DnewSAA + \sum_{\wavelength \in \waveSet}\sum_{\edge\in \edgeSetw}\dualhat^\eqref{eq: y_bound_scenarios_bd_sp_beta}_{\edge,\wavelength},$$
where $ \dualhat^{\eqref{eq: conflict_link_demands_scenarios_bd_sp_beta}}, \dualhat^{\eqref{eq: demand_UB_scenarios_bd_sp_beta_2}}, \dualhat^\eqref{eq: y_bound_scenarios_bd_sp_beta} $ are the optimal dual solutions obtained by solving the following auxiliary subproblem:
\BSE
\begin{alignat}{3}
\overline{\sub}(\hat{\beta},\xival)\ :\ \max \ \ & \sum_{(s,d) \in \SDsceSAA} \numgrantedsceSAA\label{eq: OF_bd_sp_beta_2}\\
\text{s.t.} \ \ &\sum_{\wavelength \in \waveSet}\sum_{(s,d) \in \SDsceSAA} \linksceSAA  \leq (|\edgeSetw| - \hat{\beta_\edge})  &\qquad& \edge \in \edgeSet \label{eq: conflict_link_demands_scenarios_bd_sp_beta}\\
& \eqref{eq: flow_conservation_scenarios_bd_sp_beta} - \eqref{eq: loop_elimination_scenarios_bd_sp_beta} \\
&\sum_{\wavelength \in \waveSet}\sum_{\stackrel{\edge \in}{\nodefunc^{+}(s)\cap \edgeSetw}} \linksceSAA \leq \DnewSAA &&(s,d) \in \SDsceSAA \label{eq: demand_UB_scenarios_bd_sp_beta_2} \\
&\sum_{(s,d) \in \SDsceSAA} \linksceSAA  \leq 1  && \wavelength \in \waveSet,  \edge \in \edgeSetw \label{eq: y_bound_scenarios_bd_sp_beta}\\
& \eqref{eq: granted_scenarios_sub}-\eqref{eq:granted_sub_domain}.
\end{alignat}
\ESE
Note that, $\overline{\sub}(\hat{\beta},\xival)$ is a relaxation of $\sub(\hat{x},\xival)$ where constraints \eqref{eq: conflict_link_demands_scenarios_bd_sp_beta} are obtained by integrating over constraints \eqref{eq: conflict_link_demands_scenarios_bd_sp_x_beta}. They do not prevent the wavelength conflict, hence $ \overline{\master}(\xisetSAA) $  is only a relaxation of  $ \master(\xisetSAA) $ and a complete set of $ \beta- $cuts does not give an exact reformulation of  \eqref{eq: saa-SmaxRWA}. In the next section, we  evaluate the performance of $ \beta- $cuts and show that they are very effective in reducing the number of $ x- $cuts that the algorithm needs to converge to optimality.

\section{Numerical Results}
\label{sec:num-results}
In this section, we study the effectiveness of the proposed models and algorithms for dealing with uncertainty that arises in a WDM network with dynamic traffic. In Section \ref{sec:net} we present the characteristics of the studied networks and describe the data sets. Next, in Sections \ref{sec:SmaxRWA-analysis} and \ref{sec:lr-saa}, we divide our analysis into two parts, \srwa{} and \lr{}, respectively. More specifically, for each problem we study the algorithmic performance of the decomposition framework given in Section \ref{sec:benders}, and through simulation, we compare the stochastic models with their traditional deterministic counterparts. 

\textbf{Implementation Details}. Programs are written in Java and run on Niagara\footnote{\url{https://docs.scinet.utoronto.ca/index.php/Main_Page}} servers  \citep{ponce2019deploying,loken2010scinet} using CPLEX (version V12.10.0) as the integer linear programming solver. In all experiments, a time limit of 600 seconds per \IPLP{}  is imposed on the solver. In the implementation of the decomposition algorithms, we have benefited from the $ \beta- $cuts in a pre-processing phase with first-stage decisions relaxed (\LPLP{}). The cuts are added to the \LPLP{} problem until no other $ \beta- $cut is violated. Then the first-stage decisions are converted  to integers and the \IPLP{} problem, together with the generated $ \beta- $cuts are passed to CPLEX. In this phase, we have used callbacks for adding the cuts, where $ \beta- $cuts are generated at both fractional and integral nodes, while $ x- $cuts are added only if no violated $ \beta- $cut is found at an integral node.

\subsection{Data Sets}
\label{sec:net}
\begin{figure}[h]
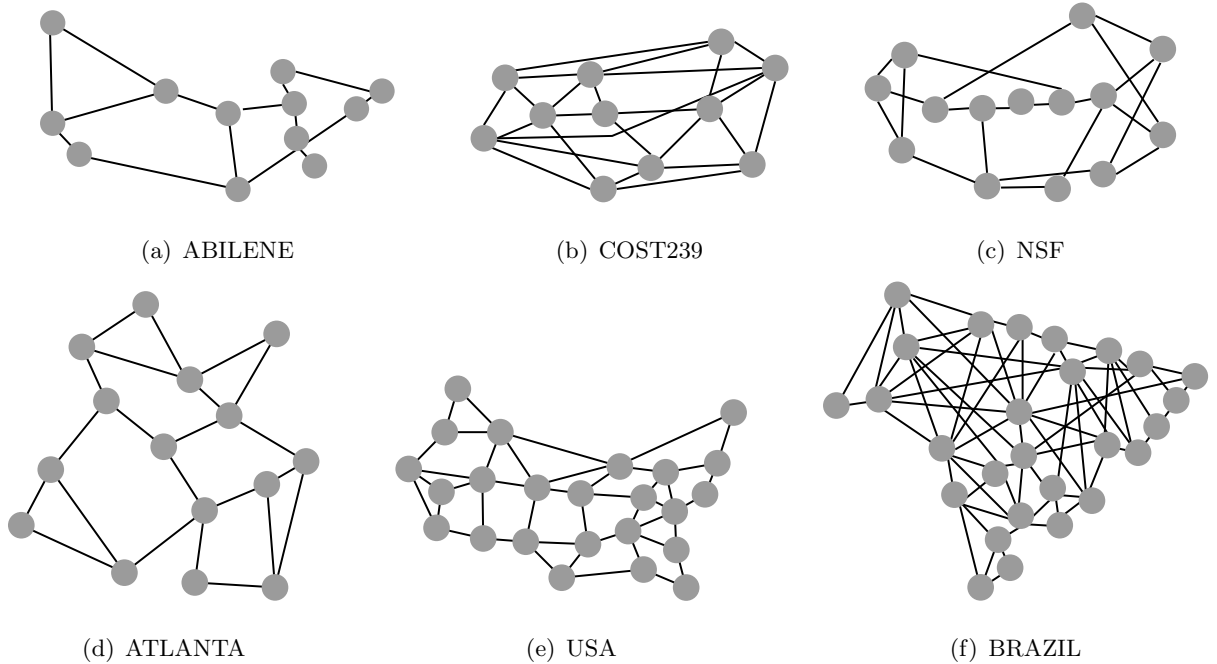

	\centering
		\subfigure[\abilene{}]{\tikzset{every picture/.style={line width=0.75pt}} 

\begin{tikzpicture}[scale=\textwidth/20cm,x=0.75pt,y=0.75pt,yscale=-1,xscale=1]

\draw [line width=0.75]    (234.13,109.65) -- (235.48,166.44) ;
\draw [line width=0.75]    (235.84,107.4) -- (305,149.2) ;
\draw [line width=0.75]    (235.8,168.8) -- (305,149.2) ;
\draw [line width=0.75]    (237.63,172.51) -- (252.01,187.79) ;
\draw [line width=0.75]    (258.54,189.95) -- (349.03,209.27) ;
\draw [line width=0.75]    (305,149.2) -- (341,161.6) ;
\draw [line width=0.75]    (343,162.8) -- (349.03,209.27) ;
\draw [line width=0.75]    (352,160.29) -- (383.4,156.8) ;
\draw [line width=0.75]    (349.03,209.27) -- (386.2,184) ;
\draw [line width=0.75]    (386.2,184) -- (391.4,189.6) ;
\draw [line width=0.75]    (383.4,156.8) -- (386.2,184) ;
\draw [line width=0.75]    (376.51,144.95) -- (383.4,156.8) ;
\draw [line width=0.75]    (386.2,184) -- (421.4,160) ;
\draw [line width=0.75]    (384.23,137.23) -- (433.8,150.4) ;
\draw [line width=0.75]    (421.4,160) -- (433.8,150.4) ;
\draw  [draw opacity=0][fill={rgb, 255:red, 155; green, 155; blue, 155 }  ,fill opacity=1 ][line width=0.75]  (228.12,107.4) .. controls (228.12,103.14) and (231.57,99.69) .. (235.84,99.69) .. controls (240.1,99.69) and (243.56,103.14) .. (243.56,107.4) .. controls (243.56,111.67) and (240.1,115.12) .. (235.84,115.12) .. controls (231.57,115.12) and (228.12,111.67) .. (228.12,107.4) -- cycle ;
\draw  [draw opacity=0][fill={rgb, 255:red, 155; green, 155; blue, 155 }  ,fill opacity=1 ][line width=0.75]  (335.28,162.8) .. controls (335.28,158.54) and (338.74,155.08) .. (343,155.08) .. controls (347.26,155.08) and (350.72,158.54) .. (350.72,162.8) .. controls (350.72,167.06) and (347.26,170.52) .. (343,170.52) .. controls (338.74,170.52) and (335.28,167.06) .. (335.28,162.8) -- cycle ;
\draw  [draw opacity=0][fill={rgb, 255:red, 155; green, 155; blue, 155 }  ,fill opacity=1 ][line width=0.75]  (228.08,168.8) .. controls (228.08,164.54) and (231.54,161.08) .. (235.8,161.08) .. controls (240.06,161.08) and (243.52,164.54) .. (243.52,168.8) .. controls (243.52,173.06) and (240.06,176.52) .. (235.8,176.52) .. controls (231.54,176.52) and (228.08,173.06) .. (228.08,168.8) -- cycle ;
\draw  [draw opacity=0][fill={rgb, 255:red, 155; green, 155; blue, 155 }  ,fill opacity=1 ][line width=0.75]  (244.3,187.79) .. controls (244.3,183.52) and (247.75,180.07) .. (252.01,180.07) .. controls (256.28,180.07) and (259.73,183.52) .. (259.73,187.79) .. controls (259.73,192.05) and (256.28,195.51) .. (252.01,195.51) .. controls (247.75,195.51) and (244.3,192.05) .. (244.3,187.79) -- cycle ;
\draw  [draw opacity=0][fill={rgb, 255:red, 155; green, 155; blue, 155 }  ,fill opacity=1 ][line width=0.75]  (341.31,209.27) .. controls (341.31,205) and (344.76,201.55) .. (349.03,201.55) .. controls (353.29,201.55) and (356.74,205) .. (356.74,209.27) .. controls (356.74,213.53) and (353.29,216.99) .. (349.03,216.99) .. controls (344.76,216.99) and (341.31,213.53) .. (341.31,209.27) -- cycle ;
\draw  [draw opacity=0][fill={rgb, 255:red, 155; green, 155; blue, 155 }  ,fill opacity=1 ][line width=0.75]  (297.28,149.2) .. controls (297.28,144.94) and (300.74,141.48) .. (305,141.48) .. controls (309.26,141.48) and (312.72,144.94) .. (312.72,149.2) .. controls (312.72,153.46) and (309.26,156.92) .. (305,156.92) .. controls (300.74,156.92) and (297.28,153.46) .. (297.28,149.2) -- cycle ;
\draw  [draw opacity=0][fill={rgb, 255:red, 155; green, 155; blue, 155 }  ,fill opacity=1 ][line width=0.75]  (375.68,156.8) .. controls (375.68,152.54) and (379.14,149.08) .. (383.4,149.08) .. controls (387.66,149.08) and (391.12,152.54) .. (391.12,156.8) .. controls (391.12,161.06) and (387.66,164.52) .. (383.4,164.52) .. controls (379.14,164.52) and (375.68,161.06) .. (375.68,156.8) -- cycle ;
\draw  [draw opacity=0][fill={rgb, 255:red, 155; green, 155; blue, 155 }  ,fill opacity=1 ][line width=0.75]  (388.07,194.9) .. controls (388.07,190.64) and (391.53,187.19) .. (395.79,187.19) .. controls (400.05,187.19) and (403.51,190.64) .. (403.51,194.9) .. controls (403.51,199.17) and (400.05,202.62) .. (395.79,202.62) .. controls (391.53,202.62) and (388.07,199.17) .. (388.07,194.9) -- cycle ;
\draw  [draw opacity=0][fill={rgb, 255:red, 155; green, 155; blue, 155 }  ,fill opacity=1 ][line width=0.75]  (377.08,178.12) .. controls (377.08,173.86) and (380.54,170.4) .. (384.8,170.4) .. controls (389.06,170.4) and (392.52,173.86) .. (392.52,178.12) .. controls (392.52,182.38) and (389.06,185.84) .. (384.8,185.84) .. controls (380.54,185.84) and (377.08,182.38) .. (377.08,178.12) -- cycle ;
\draw  [draw opacity=0][fill={rgb, 255:red, 155; green, 155; blue, 155 }  ,fill opacity=1 ][line width=0.75]  (368.79,137.23) .. controls (368.79,132.97) and (372.25,129.51) .. (376.51,129.51) .. controls (380.77,129.51) and (384.23,132.97) .. (384.23,137.23) .. controls (384.23,141.49) and (380.77,144.95) .. (376.51,144.95) .. controls (372.25,144.95) and (368.79,141.49) .. (368.79,137.23) -- cycle ;
\draw  [draw opacity=0][fill={rgb, 255:red, 155; green, 155; blue, 155 }  ,fill opacity=1 ][line width=0.75]  (413.68,160) .. controls (413.68,155.74) and (417.14,152.28) .. (421.4,152.28) .. controls (425.66,152.28) and (429.12,155.74) .. (429.12,160) .. controls (429.12,164.26) and (425.66,167.72) .. (421.4,167.72) .. controls (417.14,167.72) and (413.68,164.26) .. (413.68,160) -- cycle ;
\draw  [draw opacity=0][fill={rgb, 255:red, 155; green, 155; blue, 155 }  ,fill opacity=1 ][line width=0.75]  (429.3,149.1) .. controls (429.3,144.84) and (432.75,141.39) .. (437.02,141.39) .. controls (441.28,141.39) and (444.73,144.84) .. (444.73,149.1) .. controls (444.73,153.37) and (441.28,156.82) .. (437.02,156.82) .. controls (432.75,156.82) and (429.3,153.37) .. (429.3,149.1) -- cycle ;

\end{tikzpicture}} $ \quad $
		\subfigure[\cost{}]{\tikzset{every picture/.style={line width=0.75pt}} 

\begin{tikzpicture}[scale=\textwidth/20cm,x=0.75pt,y=0.75pt,yscale=-1,xscale=1]

\draw [color={rgb, 255:red, 0; green, 0; blue, 0 }  ,draw opacity=1 ][line width=0.75]    (199.61,108.09) -- (186.61,128.83) ;
\draw [color={rgb, 255:red, 0; green, 0; blue, 0 }  ,draw opacity=1 ][line width=0.75]    (323.48,77.96) -- (205,95.5) ;
\draw [color={rgb, 255:red, 0; green, 0; blue, 0 }  ,draw opacity=1 ][line width=0.75]    (218,117.5) -- (205,105.5) ;
\draw [color={rgb, 255:red, 0; green, 0; blue, 0 }  ,draw opacity=1 ][line width=0.75]    (243.48,97.96) -- (207.74,99.96) ;
\draw [color={rgb, 255:red, 0; green, 0; blue, 0 }  ,draw opacity=1 ][line width=0.75]    (217,127.5) -- (194.74,136.96) ;
\draw [color={rgb, 255:red, 0; green, 0; blue, 0 }  ,draw opacity=1 ][line width=0.75]    (251.48,167.96) -- (191,142.5) ;
\draw [color={rgb, 255:red, 0; green, 0; blue, 0 }  ,draw opacity=1 ][line width=0.75]    (347,156.5) -- (267.74,167.96) ;
\draw [color={rgb, 255:red, 0; green, 0; blue, 0 }  ,draw opacity=1 ][line width=0.75]    (288.61,154.96) -- (186.61,136.96) ;
\draw [color={rgb, 255:red, 0; green, 0; blue, 0 }  ,draw opacity=1 ][line width=0.75]    (288.61,154.96) -- (259.61,167.96) ;
\draw [color={rgb, 255:red, 0; green, 0; blue, 0 }  ,draw opacity=1 ][line width=0.75]    (222.61,122.96) -- (259.61,167.96) ;
\draw [color={rgb, 255:red, 0; green, 0; blue, 0 }  ,draw opacity=1 ][line width=0.75]    (251.61,97.96) -- (222.61,122.96) ;
\draw [color={rgb, 255:red, 0; green, 0; blue, 0 }  ,draw opacity=1 ][line width=0.75]    (331.61,77.96) -- (251.61,97.96) ;
\draw [color={rgb, 255:red, 0; green, 0; blue, 0 }  ,draw opacity=1 ][line width=0.75]    (356.48,93.96) -- (259.74,97.96) ;
\draw [color={rgb, 255:red, 0; green, 0; blue, 0 }  ,draw opacity=1 ][line width=0.75]    (342.48,152.96) -- (296.74,154.96) ;
\draw [color={rgb, 255:red, 0; green, 0; blue, 0 }  ,draw opacity=1 ][line width=0.75]    (260.61,121.96) -- (288.61,146.83) ;
\draw [color={rgb, 255:red, 0; green, 0; blue, 0 }  ,draw opacity=1 ][line width=0.75]    (324.61,118.96) -- (293,150.5) ;
\draw [color={rgb, 255:red, 0; green, 0; blue, 0 }  ,draw opacity=1 ][line width=0.75]    (260.61,121.96) -- (230.74,122.96) ;
\draw [color={rgb, 255:red, 0; green, 0; blue, 0 }  ,draw opacity=1 ][line width=0.75]    (260.61,121.96) -- (251.61,97.96) ;
\draw [color={rgb, 255:red, 0; green, 0; blue, 0 }  ,draw opacity=1 ][line width=0.75]    (324.61,118.96) -- (260.61,121.96) ;
\draw [color={rgb, 255:red, 0; green, 0; blue, 0 }  ,draw opacity=1 ][line width=0.75]    (324.61,110.83) -- (331.61,86.09) ;
\draw [color={rgb, 255:red, 0; green, 0; blue, 0 }  ,draw opacity=1 ][line width=0.75]    (364.61,93.96) -- (339.74,77.96) ;
\draw [color={rgb, 255:red, 0; green, 0; blue, 0 }  ,draw opacity=1 ][line width=0.75]    (350.61,152.96) -- (364.61,102.09) ;
\draw [color={rgb, 255:red, 0; green, 0; blue, 0 }  ,draw opacity=1 ][line width=0.75]    (350.61,152.96) -- (329.61,122.96) ;
\draw [color={rgb, 255:red, 0; green, 0; blue, 0 }  ,draw opacity=1 ][line width=0.75]    (324.61,118.96) -- (364.61,93.96) ;
\draw [color={rgb, 255:red, 0; green, 0; blue, 0 }  ,draw opacity=1 ][line width=0.75]    (194.74,136.96) -- (265,135.5) -- (364.61,93.96) ;
\draw  [draw opacity=0][fill={rgb, 255:red, 155; green, 155; blue, 155 }  ,fill opacity=1 ][line width=0.75]  (191.48,99.96) .. controls (191.48,95.47) and (195.12,91.83) .. (199.61,91.83) .. controls (204.1,91.83) and (207.74,95.47) .. (207.74,99.96) .. controls (207.74,104.45) and (204.1,108.09) .. (199.61,108.09) .. controls (195.12,108.09) and (191.48,104.45) .. (191.48,99.96) -- cycle ;
\draw  [draw opacity=0][fill={rgb, 255:red, 155; green, 155; blue, 155 }  ,fill opacity=1 ][line width=0.75]  (178.48,136.96) .. controls (178.48,132.47) and (182.12,128.83) .. (186.61,128.83) .. controls (191.1,128.83) and (194.74,132.47) .. (194.74,136.96) .. controls (194.74,141.45) and (191.1,145.09) .. (186.61,145.09) .. controls (182.12,145.09) and (178.48,141.45) .. (178.48,136.96) -- cycle ;
\draw  [draw opacity=0][fill={rgb, 255:red, 155; green, 155; blue, 155 }  ,fill opacity=1 ][line width=0.75]  (214.48,122.96) .. controls (214.48,118.47) and (218.12,114.83) .. (222.61,114.83) .. controls (227.1,114.83) and (230.74,118.47) .. (230.74,122.96) .. controls (230.74,127.45) and (227.1,131.09) .. (222.61,131.09) .. controls (218.12,131.09) and (214.48,127.45) .. (214.48,122.96) -- cycle ;
\draw  [draw opacity=0][fill={rgb, 255:red, 155; green, 155; blue, 155 }  ,fill opacity=1 ][line width=0.75]  (243.48,97.96) .. controls (243.48,93.47) and (247.12,89.83) .. (251.61,89.83) .. controls (256.1,89.83) and (259.74,93.47) .. (259.74,97.96) .. controls (259.74,102.45) and (256.1,106.09) .. (251.61,106.09) .. controls (247.12,106.09) and (243.48,102.45) .. (243.48,97.96) -- cycle ;
\draw  [draw opacity=0][fill={rgb, 255:red, 155; green, 155; blue, 155 }  ,fill opacity=1 ][line width=0.75]  (252.48,121.96) .. controls (252.48,117.47) and (256.12,113.83) .. (260.61,113.83) .. controls (265.1,113.83) and (268.74,117.47) .. (268.74,121.96) .. controls (268.74,126.45) and (265.1,130.09) .. (260.61,130.09) .. controls (256.12,130.09) and (252.48,126.45) .. (252.48,121.96) -- cycle ;
\draw  [draw opacity=0][fill={rgb, 255:red, 155; green, 155; blue, 155 }  ,fill opacity=1 ][line width=0.75]  (323.48,77.96) .. controls (323.48,73.47) and (327.12,69.83) .. (331.61,69.83) .. controls (336.1,69.83) and (339.74,73.47) .. (339.74,77.96) .. controls (339.74,82.45) and (336.1,86.09) .. (331.61,86.09) .. controls (327.12,86.09) and (323.48,82.45) .. (323.48,77.96) -- cycle ;
\draw  [draw opacity=0][fill={rgb, 255:red, 155; green, 155; blue, 155 }  ,fill opacity=1 ][line width=0.75]  (356.48,93.96) .. controls (356.48,89.47) and (360.12,85.83) .. (364.61,85.83) .. controls (369.1,85.83) and (372.74,89.47) .. (372.74,93.96) .. controls (372.74,98.45) and (369.1,102.09) .. (364.61,102.09) .. controls (360.12,102.09) and (356.48,98.45) .. (356.48,93.96) -- cycle ;
\draw  [draw opacity=0][fill={rgb, 255:red, 155; green, 155; blue, 155 }  ,fill opacity=1 ][line width=0.75]  (316.48,118.96) .. controls (316.48,114.47) and (320.12,110.83) .. (324.61,110.83) .. controls (329.1,110.83) and (332.74,114.47) .. (332.74,118.96) .. controls (332.74,123.45) and (329.1,127.09) .. (324.61,127.09) .. controls (320.12,127.09) and (316.48,123.45) .. (316.48,118.96) -- cycle ;
\draw  [draw opacity=0][fill={rgb, 255:red, 155; green, 155; blue, 155 }  ,fill opacity=1 ][line width=0.75]  (342.48,152.96) .. controls (342.48,148.47) and (346.12,144.83) .. (350.61,144.83) .. controls (355.1,144.83) and (358.74,148.47) .. (358.74,152.96) .. controls (358.74,157.45) and (355.1,161.09) .. (350.61,161.09) .. controls (346.12,161.09) and (342.48,157.45) .. (342.48,152.96) -- cycle ;
\draw  [draw opacity=0][fill={rgb, 255:red, 155; green, 155; blue, 155 }  ,fill opacity=1 ][line width=0.75]  (280.48,154.96) .. controls (280.48,150.47) and (284.12,146.83) .. (288.61,146.83) .. controls (293.1,146.83) and (296.74,150.47) .. (296.74,154.96) .. controls (296.74,159.45) and (293.1,163.09) .. (288.61,163.09) .. controls (284.12,163.09) and (280.48,159.45) .. (280.48,154.96) -- cycle ;
\draw  [draw opacity=0][fill={rgb, 255:red, 155; green, 155; blue, 155 }  ,fill opacity=1 ][line width=0.75]  (251.48,167.96) .. controls (251.48,163.47) and (255.12,159.83) .. (259.61,159.83) .. controls (264.1,159.83) and (267.74,163.47) .. (267.74,167.96) .. controls (267.74,172.45) and (264.1,176.09) .. (259.61,176.09) .. controls (255.12,176.09) and (251.48,172.45) .. (251.48,167.96) -- cycle ;

\end{tikzpicture}} $ \quad $ 
		\subfigure[\nsf{}]{\tikzset{every picture/.style={line width=0.75pt}} 

\begin{tikzpicture}[scale=\textwidth/20cm,x=0.75pt,y=0.75pt,yscale=-1,xscale=1]

\draw [line width=0.75]    (175.61,70.83) -- (184.15,61.07) ;
\draw [line width=0.75]    (182.93,81.2) -- (203.67,87.7) ;
\draw [line width=0.75]    (175.61,87.09) -- (188.22,109.87) ;
\draw [line width=0.75]    (191.88,66.76) -- (190.25,108.65) ;
\draw [line width=0.75]    (218.72,91.97) -- (231.32,91.16) ;
\draw [line width=0.75]    (200.01,58.63) -- (287.85,79.37) ;
\draw  [draw opacity=0][fill={rgb, 255:red, 155; green, 155; blue, 155 }  ,fill opacity=1 ][line width=0.75]  (167.48,78.96) .. controls (167.48,74.47) and (171.12,70.83) .. (175.61,70.83) .. controls (180.1,70.83) and (183.74,74.47) .. (183.74,78.96) .. controls (183.74,83.45) and (180.1,87.09) .. (175.61,87.09) .. controls (171.12,87.09) and (167.48,83.45) .. (167.48,78.96) -- cycle ;
\draw  [draw opacity=0][fill={rgb, 255:red, 155; green, 155; blue, 155 }  ,fill opacity=1 ][line width=0.75]  (183.74,58.63) .. controls (183.74,54.13) and (187.39,50.49) .. (191.88,50.49) .. controls (196.37,50.49) and (200.01,54.13) .. (200.01,58.63) .. controls (200.01,63.12) and (196.37,66.76) .. (191.88,66.76) .. controls (187.39,66.76) and (183.74,63.12) .. (183.74,58.63) -- cycle ;
\draw  [draw opacity=0][fill={rgb, 255:red, 155; green, 155; blue, 155 }  ,fill opacity=1 ][line width=0.75]  (182.12,116.78) .. controls (182.12,112.29) and (185.76,108.65) .. (190.25,108.65) .. controls (194.74,108.65) and (198.38,112.29) .. (198.38,116.78) .. controls (198.38,121.27) and (194.74,124.91) .. (190.25,124.91) .. controls (185.76,124.91) and (182.12,121.27) .. (182.12,116.78) -- cycle ;
\draw  [draw opacity=0][fill={rgb, 255:red, 155; green, 155; blue, 155 }  ,fill opacity=1 ][line width=0.75]  (202.45,91.97) .. controls (202.45,87.48) and (206.09,83.84) .. (210.58,83.84) .. controls (215.08,83.84) and (218.72,87.48) .. (218.72,91.97) .. controls (218.72,96.47) and (215.08,100.11) .. (210.58,100.11) .. controls (206.09,100.11) and (202.45,96.47) .. (202.45,91.97) -- cycle ;
\draw  [draw opacity=0][fill={rgb, 255:red, 155; green, 155; blue, 155 }  ,fill opacity=1 ][line width=0.75]  (231.32,91.16) .. controls (231.32,86.67) and (234.97,83.03) .. (239.46,83.03) .. controls (243.95,83.03) and (247.59,86.67) .. (247.59,91.16) .. controls (247.59,95.65) and (243.95,99.29) .. (239.46,99.29) .. controls (234.97,99.29) and (231.32,95.65) .. (231.32,91.16) -- cycle ;
\draw  [draw opacity=0][fill={rgb, 255:red, 155; green, 155; blue, 155 }  ,fill opacity=1 ][line width=0.75]  (234.17,138.74) .. controls (234.17,134.25) and (237.81,130.61) .. (242.3,130.61) .. controls (246.8,130.61) and (250.44,134.25) .. (250.44,138.74) .. controls (250.44,143.23) and (246.8,146.87) .. (242.3,146.87) .. controls (237.81,146.87) and (234.17,143.23) .. (234.17,138.74) -- cycle ;
\draw  [draw opacity=0][fill={rgb, 255:red, 155; green, 155; blue, 155 }  ,fill opacity=1 ][line width=0.75]  (254.91,87.09) .. controls (254.91,82.6) and (258.55,78.96) .. (263.04,78.96) .. controls (267.54,78.96) and (271.18,82.6) .. (271.18,87.09) .. controls (271.18,91.59) and (267.54,95.23) .. (263.04,95.23) .. controls (258.55,95.23) and (254.91,91.59) .. (254.91,87.09) -- cycle ;
\draw  [draw opacity=0][fill={rgb, 255:red, 155; green, 155; blue, 155 }  ,fill opacity=1 ][line width=0.75]  (279.72,87.5) .. controls (279.72,83.01) and (283.36,79.37) .. (287.85,79.37) .. controls (292.34,79.37) and (295.98,83.01) .. (295.98,87.5) .. controls (295.98,91.99) and (292.34,95.63) .. (287.85,95.63) .. controls (283.36,95.63) and (279.72,91.99) .. (279.72,87.5) -- cycle ;
\draw  [draw opacity=0][fill={rgb, 255:red, 155; green, 155; blue, 155 }  ,fill opacity=1 ][line width=0.75]  (292.32,34.63) .. controls (292.32,30.14) and (295.97,26.5) .. (300.46,26.5) .. controls (304.95,26.5) and (308.59,30.14) .. (308.59,34.63) .. controls (308.59,39.13) and (304.95,42.77) .. (300.46,42.77) .. controls (295.97,42.77) and (292.32,39.13) .. (292.32,34.63) -- cycle ;
\draw  [draw opacity=0][fill={rgb, 255:red, 155; green, 155; blue, 155 }  ,fill opacity=1 ][line width=0.75]  (305.34,83.43) .. controls (305.34,78.94) and (308.98,75.3) .. (313.47,75.3) .. controls (317.96,75.3) and (321.6,78.94) .. (321.6,83.43) .. controls (321.6,87.93) and (317.96,91.57) .. (313.47,91.57) .. controls (308.98,91.57) and (305.34,87.93) .. (305.34,83.43) -- cycle ;
\draw  [draw opacity=0][fill={rgb, 255:red, 155; green, 155; blue, 155 }  ,fill opacity=1 ][line width=0.75]  (277.28,140.37) .. controls (277.28,135.87) and (280.92,132.23) .. (285.41,132.23) .. controls (289.9,132.23) and (293.54,135.87) .. (293.54,140.37) .. controls (293.54,144.86) and (289.9,148.5) .. (285.41,148.5) .. controls (280.92,148.5) and (277.28,144.86) .. (277.28,140.37) -- cycle ;
\draw  [draw opacity=0][fill={rgb, 255:red, 155; green, 155; blue, 155 }  ,fill opacity=1 ][line width=0.75]  (341.53,54.97) .. controls (341.53,50.47) and (345.17,46.83) .. (349.66,46.83) .. controls (354.16,46.83) and (357.8,50.47) .. (357.8,54.97) .. controls (357.8,59.46) and (354.16,63.1) .. (349.66,63.1) .. controls (345.17,63.1) and (341.53,59.46) .. (341.53,54.97) -- cycle ;
\draw  [draw opacity=0][fill={rgb, 255:red, 155; green, 155; blue, 155 }  ,fill opacity=1 ][line width=0.75]  (304.93,131.42) .. controls (304.93,126.93) and (308.57,123.29) .. (313.06,123.29) .. controls (317.56,123.29) and (321.2,126.93) .. (321.2,131.42) .. controls (321.2,135.91) and (317.56,139.55) .. (313.06,139.55) .. controls (308.57,139.55) and (304.93,135.91) .. (304.93,131.42) -- cycle ;
\draw  [draw opacity=0][fill={rgb, 255:red, 155; green, 155; blue, 155 }  ,fill opacity=1 ][line width=0.75]  (341.53,107.43) .. controls (341.53,102.93) and (345.17,99.29) .. (349.66,99.29) .. controls (354.16,99.29) and (357.8,102.93) .. (357.8,107.43) .. controls (357.8,111.92) and (354.16,115.56) .. (349.66,115.56) .. controls (345.17,115.56) and (341.53,111.92) .. (341.53,107.43) -- cycle ;

\draw [line width=0.75]    (247.59,91.16) -- (256.13,91.06) ;
\draw [line width=0.75]    (270.97,89.03) -- (280.12,88.62) ;
\draw [line width=0.75]    (295.98,87.5) -- (305.74,85.57) ;
\draw [line width=0.75]    (320.99,80.28) -- (344.17,61.37) ;
\draw [line width=0.75]    (319.37,89.03) -- (342.14,105.5) ;
\draw [line width=0.75]    (215.87,85.26) -- (293.54,39.21) ;
\draw [line width=0.75]    (308.59,34.63) -- (342.75,50.6) ;
\draw [line width=0.75]    (306.56,39.61) -- (346.82,99.6) ;
\draw [line width=0.75]    (317.74,124.61) -- (348.24,63.41) ;
\draw [line width=0.75]    (312.25,91.47) -- (289.27,132.74) ;
\draw [line width=0.75]    (239.46,99.29) -- (242.3,130.61) ;
\draw [line width=0.75]    (196.96,120.74) -- (234.98,136.4) ;
\draw [line width=0.75]    (277.48,139.25) -- (250.84,139.65) ;
\draw [line width=0.75]    (305.34,129.29) -- (249.83,134.98) ;
\draw [line width=0.75]    (321.2,131.42) -- (347.22,115.36) ;

\end{tikzpicture}}\\
		\subfigure[\atlanta{}]{\input{./Figures/atlanta.tex}} $ \quad $
		\subfigure[\usa{}]{\input{./Figures/USA.tex}} $ \quad $
		\subfigure[\brazil{}]{\input{./Figures/brazil.tex}} 
		\caption{WDM Networks}
		\label{fig:networks}
\end{figure}

For the physical topology of our instances, we have considered six standard long-haul networks (Table \ref{tab:networks}) from the literature of WDM networks.  These topologies represent different levels of connectivity degree (column ``avg deg.''), ranging  from \abilene{} with a small number of links and paths, to the highly interconnected network \brazil{}. In Figure \ref{fig:networks} we have reproduced the networks to give a visual representation of their meshness. 

\begin{table}[htbp]
	\small
  \centering
  \caption{Networks}
\begin{tabular}{lrrrrl}
	\toprule
Network & \multicolumn{1}{l}{$|\nodeSet|$} & \multicolumn{1}{l}{$|\edgeSet|$} & \multicolumn{1}{l}{avg deg.} & Source \\
\midrule
\abilene{} & 12    & 30    & 2.5      & \cite{orlowski2010sndlib} \\
\cost{} & 11    & 50    & 4.5   &\cite{tan1996wavelength} \\
\nsf{}   & 14    & 42    & 3.0      & \cite{ramaswami1996design} \\
\atlanta{} & 15    & 44    & 2.9      & \cite{orlowski2010sndlib} \\
\usa{}   & 24    & 88    & 3.7    & \cite{batayneh2011routing} \\
\brazil{} & 27    & 140   & 5.2          & \cite{noronha2006routing} \\
\bottomrule
\end{tabular}%
  \label{tab:networks}%
\end{table}%

Arrivals of the requests (size of the incoming batches) follow a Poisson distribution with parameter $ \arrivals $, while for service holding times (how long a connection remains on the network), we have considered an exponential distribution with mean $ 1/\drops $ (as in \cite{chen2015spectrum} and \cite{xiong2018sdn}). Distribution of lightpath requests among node pairs is uniform (the same as \cite{giorgetti2015dynamic}). For each network, the arrival and service holding time  parameters (provided for each analysis in their respective subsection) are selected via some preliminary experiments to avoid trivial solutions, i.e., such that the  decisions are visibly affected by the uncertainty of incoming traffic \footnote{\revacc{The generated instances, along with the detailed results on the GoS and link usage, are available for download at \url{https://github.com/mdaryalal/DataFiles/tree/master/StochasticRWA}}}. 

Note that, if the size of the batches with respect to the capacity of the network and the length of the planning horizon is too small, the assignment and routing decisions matter less, as the network would be free to accept all the requests regardless of the previous decisions. On the other hand, if the size is too large, the network becomes rapidly congested and soon there is no room to grant any new requests via deterministic or stochastic models. For this reason, for each experiment we have chosen the parameters $ \arrivals $ and $ |\waveSet| $  based on the considered problem (\srwa{} or \lr{}), network and the number of requests on the network. 

\subsection{Analysis of  the \srwa{} Problem}\label{sec:SmaxRWA-analysis}
In this section, we study the model and algorithms proposed for \srwa{} for the networks given in Table \ref{tab:networks}. 

\subsubsection{Algorithmic Performance}\label{sec:SmaxRWA-perf}
We have evaluated the performance of the decomposition method in comparison to solving the extensive form (see Section \ref{sec:saa}) with CPLEX, as well as the strength of the $ \beta-$cut family. For this analysis, the smallest network \abilene{} is selected (for other parameter settings, see Section \ref{sec:simul-SmaxRWA}). The evaluation is in terms of the solution time in seconds, optimality gap at the end of the time limit of 10 minutes, and in case of the decomposition methods, the number of added cuts. In Table \ref{tab:maxrwa-alg-perf}, the average results of 30 repetitions for each scenario level is reported. In this table BENDERS-$x$ stands for the Benders decomposition method when it purely relies on the $ x-$cuts, while BENDERS-$x\beta $ embeds both family of cuts. If an instance does not converge within the time limit, it is indicated by ``TL'', and if it does not return a gap, it is marked by ``NA''.  To make the comparison more meaningful, in this experiment, we omitted the \LPLP{} phase and cut generation at fractional solutions. The reason for the latter is that, as demonstrated by the results in Table \ref{tab:maxrwa-alg-perf}, BENDERS-$x$ adds a large number of $ x-$cuts, and increasing them by allowing cut generation at fractional solutions heavily impacts its performance, e.g., the algorithm is not able to converge even in the smallest scenario level. 

\begin{table}[h]
	\centering
	\caption{Average algorithmic performance of various methods for solving \srwa{} problem over 30 repetitions, when \abilene{} network is considered.}
	\small
\begin{tabular}{rrrrrrrrrr}
\toprule
\multicolumn{1}{c}{\multirow{2}{*}{$|\xisetSAA|$}} & \multicolumn{2}{c}{CPLEX} & \multicolumn{3}{c}{BENDERS-$x$} & \multicolumn{4}{c}{BENDERS-$x\beta$} \\
\cmidrule(lr){2-3}  \cmidrule(lr){4-6} \cmidrule(lr){7-10}    & \multicolumn{1}{c}{time (s)} & \multicolumn{1}{c}{gap (\%)} & \multicolumn{1}{c}{time (s)} & \multicolumn{1}{c}{gap (\%)} & \multicolumn{1}{c}{$\#  x-$cuts} & \multicolumn{1}{c}{time (s)} & \multicolumn{1}{c}{gap (\%)} & \multicolumn{1}{c}{$\#  \beta-$cuts} & \multicolumn{1}{c}{$\#  x-$cuts} \\
\midrule
10    & 200   & 0     & 187   & 3     & 1124  & 17    & 0     & 8     & 15 \\
20    & TL    & NA    & 465   & 2     & 1370  & 42    & 0     & 19    & 61 \\
30    & TL    & NA    & 536   & 3     & 1538  & 157   & 0     & 31    & 204 \\
40    & TL    & NA    & TL    & 3     & 1785  & 164   & 0     & 34    & 203 \\
50    & TL    & NA    & TL    & 3     & 1980  & 206   & 0     & 58    & 209 \\
100   & TL    & NA    & TL    & 3     & 2164  & 305   & 0     & 99    & 400 \\
\bottomrule
\end{tabular}%
	\label{tab:maxrwa-alg-perf}%
\end{table}%

The results of Table \ref{tab:maxrwa-alg-perf} clearly show the superiority of BENDERS-$x\beta$ compared to the other two methods. As expected, CPLEX is unable to solve the extensive form with more than a few scenarios, and after the first scenario level, it even fails to provide a  nonzero feasible solution with an optimality gap. BENDERS-$x\beta$ solves every instance in the given time limit. It is able to significantly reduce the number of $ x-$cuts by adding a small number of $ \beta-$cuts, demonstrating the strength of this new family of cuts.

\subsubsection{Comparison with the  \maxRWA{} Problem}\label{sec:simul-SmaxRWA}
\begin{figure}[htbp]
	\centering
	\small
	\input{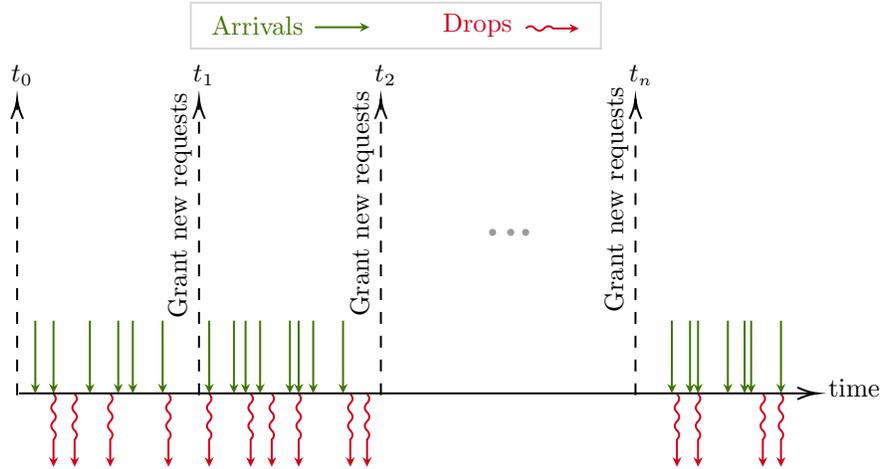}%
	\caption{Simulation framework for \srwa{}}
	\label{fig:simul-frame-SmaxRWA}
\end{figure}

Over a planning horizon of 52 stages (representing weeks of a year), we have simulated the outcomes of \srwa{} and \maxRWA{} in granting new requests on the networks of Table \ref{tab:networks}, when we have batch arrivals and predetermined decision-making times. Networks are initially empty and the simulation is repeated for 50 times. Thus, for each instance, 50 sample paths are generated, which are characterized by the incoming requests at every stage, along with their service holding time. For the latter, $ \drops\in\{13,26,52\} $ are used in the experiments for each network. At every stage, \srwa{} and \maxRWA{} grant new requests by assigning lightpaths. As the stages move forward, some requests are dropped according to their service holding times, releasing resources that the models can use. This process is outlined in Figure \ref{fig:simul-frame-SmaxRWA}.

Our decisions matter the most at certain ratios of arrival rate to the number of wavelengths, $\arrivals/|\waveSet| $, as a measure of the growth rate of the GoS, specifically in those where the network neither remains empty enough to accept all the requests over the planning horizon, nor quickly gets congested without any bandwidth room to grant new requests. In either of these two states, both models perform practically the same. By experimenting with various parameter settings, an interstate is achieved using the values in Table \ref{tab:simul-param-SmaxRWA}. In this table, we have focused on the ratio of $\arrivals/|\waveSet| $ and kept the size of $ \waveSet $ small enough such that the simulation framework is able to proceed for a long planning horizon of 52 stages. 

\begin{table}[htbp]
  \centering
  \caption{Parameter Setting for \srwa{} vs. \maxRWA{} Simulation}
  \small
\begin{tabular}{lcccccc}
\toprule
 & \abilene{} & \cost{} & \nsf{}   & \atlanta{} & \usa{}   & \brazil{} \\
\midrule
$\arrivals$ & 40    & 10    & 50    & 60    & 40    & 50 \\
$|\waveSet|$ & 10    & 10    & 10    & 10    & 5     & 5 \\
\bottomrule
\end{tabular}%
  \label{tab:simul-param-SmaxRWA}%
\end{table}%

In order to decide on the number of scenarios to be used in the SAA problems in the remaining experiments, an SAA analysis is performed for the two smallest networks, \abilene{} and \cost{}. The full SAA analysis is summarized in  
\revacc{{\App} (Section \ref{appendix:saa-results})}.
Based on the  results, sample size of $ |\xisetSAA|=500 $ is selected for \abilene{} and \cost{} instances. For the rest of the networks, the largest computationally affordable sample size in the given time limit is employed: $ |\xisetSAA|=250 $ for \nsf{}, \atlanta{} and \usa{}, $ |\xisetSAA|=50 $ for \brazil{}.

\begin{figure}[htbp]
	\centering
	\small
	\input{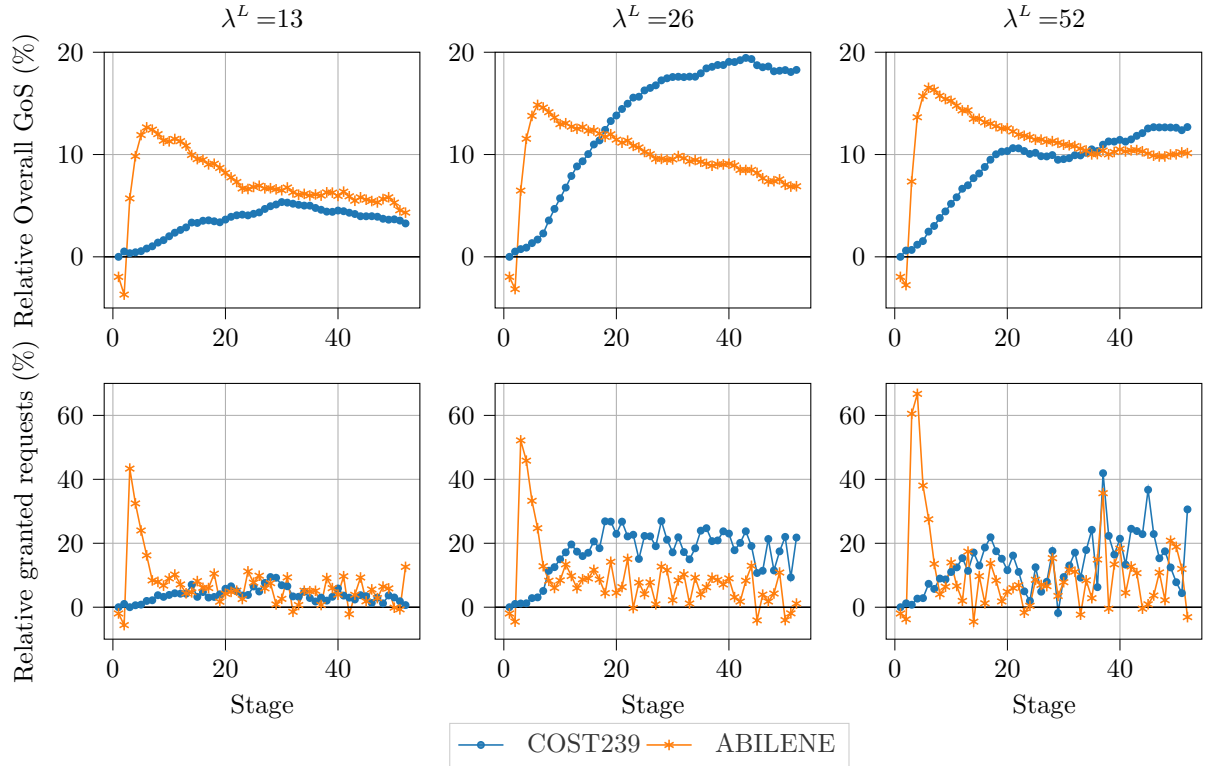}%
	\caption{Relative performance of \srwa{} for small-sized networks over different simulation horizons}
	\label{fig:improved-cost239abilenegos-add}
\end{figure}

\begin{figure}[htbp]
	\centering
	\small
	\input{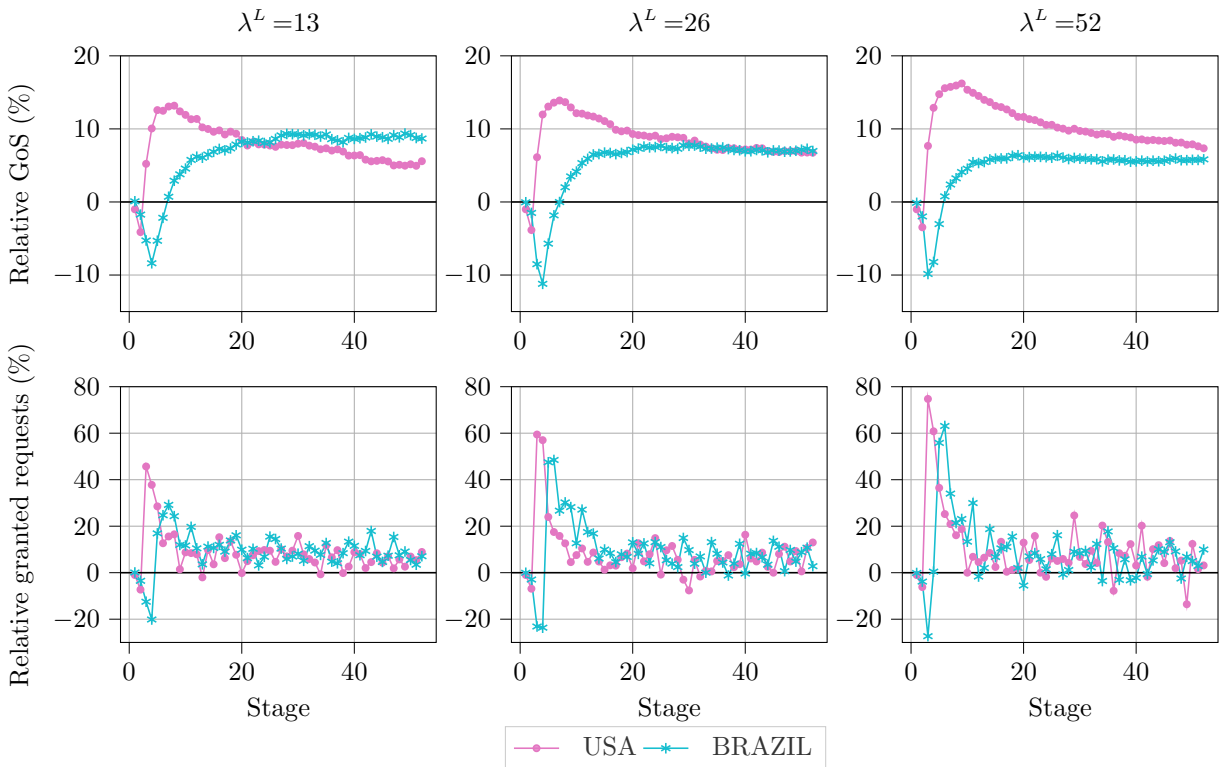}%
	\caption{Relative performance of \srwa{} for large-sized networks over different simulation horizons}
	\label{fig:improved-usabrazilgos-add}
\end{figure}

Figures \ref{fig:improved-cost239abilenegos-add} and \ref{fig:improved-usabrazilgos-add} present the relative performance of \srwa{} in terms of overall GoS and the number of granted requests at each stage, compared to the deterministic \maxRWA{} problem which does not take potential  future arrivals into consideration. These figures show the percentage of improvement or loss obtained from incorporating \srwa{} in comparison to \maxRWA{}. Results for medium-sized networks, namely \nsf{} and \atlanta{}, are given in 
\revacc{{\App} (Section \ref{appendix:simu-results})}.
In all instances, by using \srwa{} for lightpath assignment decisions, the number of granted requests and  GoS have improved over the planning horizon, with improvements up to 73\% in granted requests (\usa{} with $ \drops=52 $) and 19\% for the  GoS (\cost{} with $ \drops=26 $). \rev{The results show that  \srwa{} opts for less request grants at earlier stages to reserve the resources for future arrivals. After few stages, it catches up with and surpasses the deterministic model in terms of the GoS. Over the course of the planning horizon, for the number of granted requests we observe oscillations. It should be noted that, due to the objective differences in the deterministic and stochastic models, as well as the random request arrivals and drops, there is no expected trend and the number of granted requests can oscillate.}

In terms of granted requests, whose maximization is the actual objective function of \srwa{} and \maxRWA{}, the improvement is more significant when the service holding time of a request is longer and errors in allocation of resources by the decision maker have a longer term impact. For GoS though, \cost{} and \brazil{} exhibit dissimilar patterns. As suggested by the shape of their relative GoS curve (which, unlike the other instances, appear not to have peaked), our planning horizon does not give a full picture of their behavior.  These two networks are highly interconnected, therefore they have a large number of candidate lightpaths, i.e., high capacity for granting the requests. We conjecture that the peak of the relative GoS improvement should appear over a longer planning horizon, giving a full picture of their behavior. 
\rev{In all instances, the deterministic model is able to achieve a rejection rate close to zero for the first-stage traffic matrix. It is expected that, with the progress of the simulation, the rejection rate grows. However, our stochastic model is able to maintain a rejection rate less than 1\%. For example, for \nsf{} and \brazil{} networks the rejection rate is, on average, 0.57\% and 0.35\%, respectively.}

\revminor{
\indent As the traffic distribution might considerably vary across networks \citep{patri2020planning}, we also looked into a case where on average 20\% of the nodes are responsible for 80\% of the traffic. We call this the 20/80 distribution. 
Such instances can fall under the trivial cases as described in Section \ref{sec:net}, i.e., when the network is almost completely free or congested for the new incoming traffic, which makes cautious decision making irrelevant. The numerical experiments for \cost{} and \abilene{} with the 20/80 traffic distribution, along with the related discussions are provided in the 
\revacc{{\App} (Section \ref{appendix:simu-results})}.
The results show that for shorter service holding times, the performance of the \srwa{} and \maxRWA{} are comparable because the new incoming traffic is to some extent predictable and, with a high probability, it is simply replacing the dropped requests. For larger values of the service holding time, our decisions have long-term effects and the value of using the \srwa{}  for lightpath assignments becomes more tangible.
}

Improvements in overall GoS and granted requests are  more valuable when an important aspect of solution characteristics of both models is examined: bandwidth spectrum usage. We define the spectrum usage as the number of wavelinks utilized in granting the requests. 
From the results in Figures \ref{fig:improved-cost239abilenegos-add} and \ref{fig:improved-usabrazilgos-add}, we might wonder whether \maxRWA{} is able to catch up to \srwa{} in the long run, considering the higher number of granted requests on the network if we use \srwa{}. The answer to this question depends on the resource utilization by \srwa{}, i.e., spectrum usage. 
\begin{figure}[t]
	\centering
	\hspace*{-4mm}
	\subfigure[\abilene{}]{
		\centering
		\small
		\input{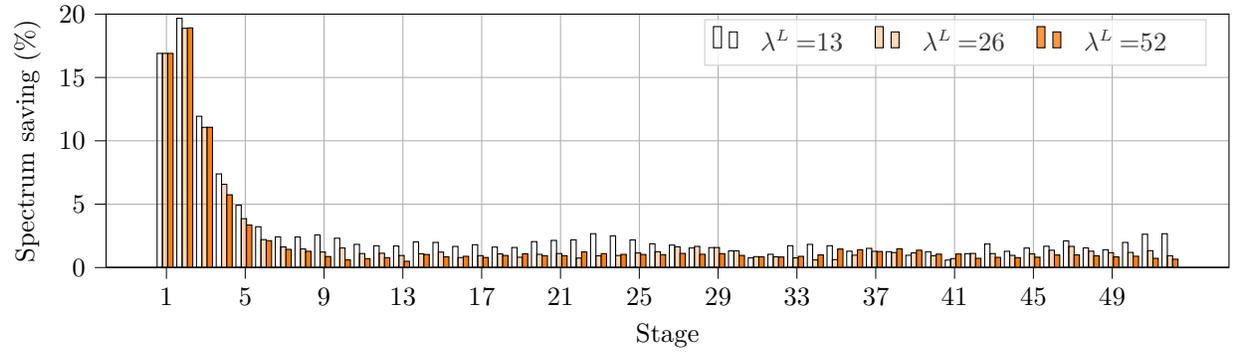}
		\label{fig:spectrum-usage-abilene}
	}

\hspace*{-4mm}
	\subfigure[\cost{}]{
		\centering
		\small
		\input{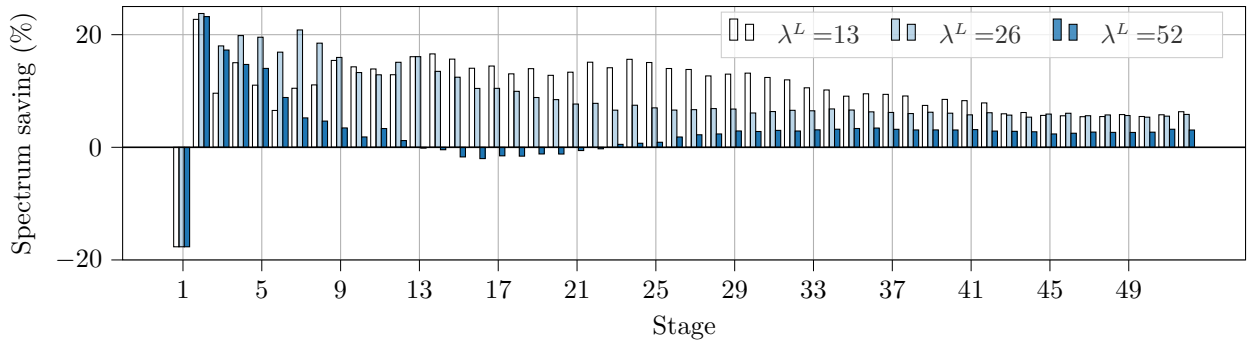}
		\label{fig:spectrum-usage-cost}
	}

\hspace*{-4mm}
	\subfigure[\brazil{}]{
		\centering
		\small
		\input{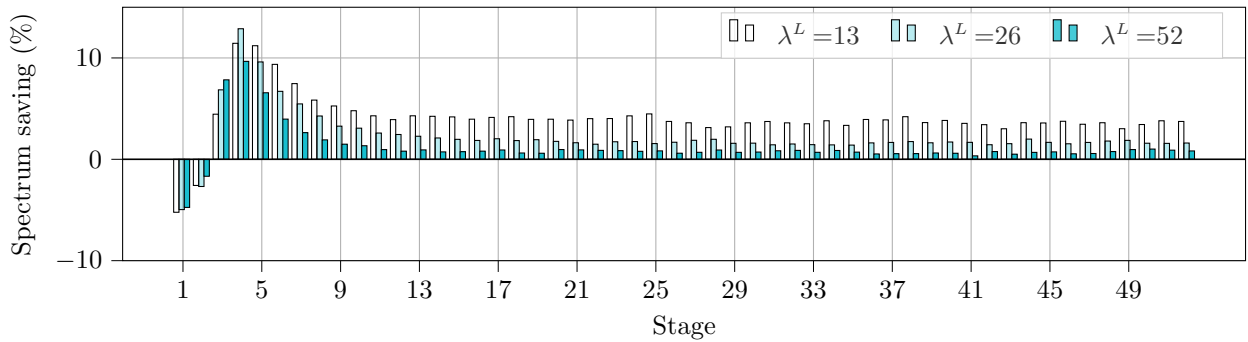}
		\label{fig:spectrum-usage-brazil}
	}

\caption{Relative bandwidth spectrum usage}
\label{fig:spec-usage}
\end{figure}
Figure \ref{fig:spec-usage} reveals that, despite the higher number of granted requests by \srwa{}, it has a lower spectrum usage, with savings of up to 20\% for the smaller networks in the earlier stages. In other words, \srwa{} achieves higher GoS by refusing some requests, as well as making informed path assignments that utilize less resources. \srwa{} effectively decreases the number of used wavelinks, an objective sought out by the \minRWA{} problem. By doing so, it postpones the fragmentation of the network and maintains a lower level of resource consumption. 
\revacc{{\App} (Section \ref{appendix:simu-results})} 
includes the set of results for the three other networks not presented in Figure \ref{fig:spec-usage}, where the same trend emerges in all of our experiments. 

\subsubsection{Quality of Relaxation}
To assess the quality of \IPLP{} bound, we have considered the network \cost{} that has the smallest $ \lambda^R $ among the instances and CPLEX is able to solve its \IPIP{} extensive form with up 50 scenarios in reasonable time. Consider the simulation framework of Figure \ref{fig:simul-frame-SmaxRWA}. For values of $ \lambda^L\in\{13,26,52\} $, we compare the optimal objective values of \IPIP{} and \IPLP{} models at different points over the planning horizon. Let $ \hat{t} $ be this particular stage. For providing a comparable setting, at each experiment, up to stage $ \hat{t} $ we grant the incoming traffic with \maxRWA{}. At stage $ \hat{t} $, we solve both \IPIP{} and \IPLP{} problems with the same sample $ \xisetSAA $ of size 50, and repeat this 30 times. Table \ref{tab:bounds-SmaxRWA} presents the comparison results. The column labeled ``$ \neq $'' gives the number of samples (out of 30) with different \IPIP{} and \IPLP{} optimal values. For these particular samples, the column ``gap'' measures the average relative gap between the two optimal values. The results show that in all cases the bound obtained by \IPLP{} problem is either exact or very close to the \IPIP{} optimal value. It is worth noting that, in every case, the first-stage cost (not necessarily the solution itself) of both models is the same, even if their optimal objective values slightly differ.  
\begin{table}[htbp]
  \centering
  \caption{Comparison of \IPLP{} bound and the optimal value of  \IPIP{}, over 30 repetitions}
    \small
    \begin{tabular}{crrr}
    \toprule
    \multicolumn{1}{c}{$ \lambda^L $} & \multicolumn{1}{c}{$\hat{t}$} & \multicolumn{1}{c}{$ \neq $} & \multicolumn{1}{c}{gap (\%)} \\
    \midrule
    \multirow{4}[2]{*}{13} & 1     & 0/30    &  NA \\
          & 13    & 1/30    & 0.11 \\
          & 26    & 0/30    &  NA\\
          & 52    & 0/30    &  NA\\
    \midrule
    \multirow{3}{*}{26} & 13    & 0/30    &  NA\\
          & 26    & 4/30    & 0.09 \\
          & 52    & 5/30    & 0.05 \\
    \midrule
    \multirow{3}[2]{*}{52} & 13    & 1/30    & 0.10 \\
          & 26    & 3/30    & 0.09 \\
          & 52    & 2/30    & 0.08 \\
    \bottomrule
    \end{tabular}%
  \label{tab:bounds-SmaxRWA}%
\end{table}%

\subsection{Analysis of the \lr{} Problem}\label{sec:lr-saa}
As a defragmentation method, \lr{} is suitable for situations where there is a large number of requests on the network, but it is not at full capacity and stranded bandwidth cannot be used in granting new requests. As such, \lr{} involves rerouting decisions for a large number of existing requests on the network, making it computationally more expensive than \srwa{}.

\subsubsection{Algorithmic Performance}
Table \ref{tab:defrag-alg-perf} presents the performance comparison of CPLEX, BENDERS-$x$ and BENDERS-$x\beta $ (see Section \ref{sec:SmaxRWA-perf}) in solving instances of \abilene{} with 100 existing requests on the network, $ |\waveSet|=10 $ and $ \arrivals=10 $. When they converge, CPLEX and BENDERS-$x$ take on the average 2.5 and 2.3 times longer than BENDERS-$x\beta $, respectively. Unlike \srwa{}, when BENDERS-$x$ does not converge within the time limit, it yields a large optimality gap, making BENDERS-$x\beta $ the sole option for solving \lr{}.
\begin{table}[h]
	\centering
	\caption{Average algorithmic performance of various methods for solving \lr{} problem over 30 repetitions, when \abilene{} network is considered.}
	\small
\begin{tabular}{rrrrrrrrrr}
\toprule
\multicolumn{1}{c}{\multirow{2}{*}{$|\xisetSAA|$}} & \multicolumn{2}{c}{CPLEX} & \multicolumn{3}{c}{BENDERS-$x$} & \multicolumn{4}{c}{BENDERS-$x\beta $} \\
\cmidrule(lr){2-3}\cmidrule(lr){4-6}\cmidrule(lr){7-10}        & \multicolumn{1}{c}{time (s)} & \multicolumn{1}{c}{gap (\%)} & \multicolumn{1}{c}{time (s)} & \multicolumn{1}{c}{gap (\%)} & \multicolumn{1}{c}{$\#  x-$cuts} & \multicolumn{1}{c}{time (s)} & \multicolumn{1}{c}{gap (\%)} & \multicolumn{1}{c}{$\# \beta-$cuts} & \multicolumn{1}{c}{$\#  x-$cuts} \\
\midrule
10    & 11    & 0     & 8     & 0     & 336   & 13    & 0     & 2     & 101 \\
20    & 120   & 0     & 82    & 0     & 847   & 67    & 0     & 7     & 426 \\
30    & 333   & 0     & 384   & 0     & 1545  & 91    & 0     & 6     & 418 \\
40    & 528   & 0     & 439   & 0     & 4315  & 139   & 0     & 7     & 877 \\
50    & TL    & NA    & TL    & 36    & 4541  & 227   & 0     & 13    & 1388 \\
100   & TL    & NA    & TL    & 57    & 7286  & 338   & 0     & 21    & 1887 \\
\bottomrule
\end{tabular}%
	\label{tab:defrag-alg-perf}%
\end{table}%

\subsubsection{Comparison with the  \minRWA{} Problem} 
\begin{figure}[t]
	\centering
	\small
	\input{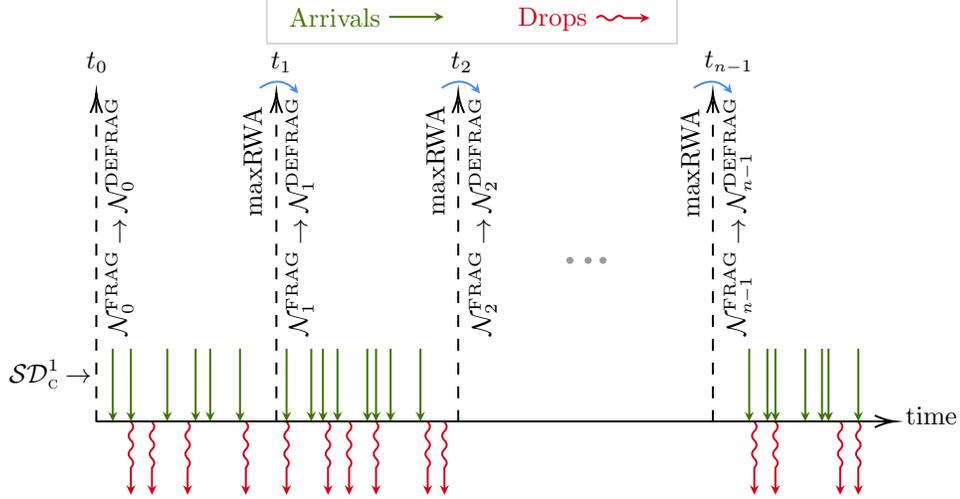}%
	\caption{Simulation framework for \lr{}}
	\label{fig:simul-frame-SmaxLR}
\end{figure}

In this section, we evaluate the impact of using \lr{} for the defragmentation of small networks \abilene{} and \cost{}. We simulate a planning horizon of 10 stages for 50 times, starting with a set of initial requests on each network. Simulation framework is depicted in Figure \ref{fig:simul-frame-SmaxLR}, where $ \mathcal{N}_t^\text{FRAG} $ and $ \mathcal{N}_t^\text{DEFRAG} $ are respectively the fragmented and defragmented provisioning at stage $ t $. Over the course of $ n $ stages, at the beginning of each stage $ t\in\{0,\dots,n-1,n\} $, a defragmentation event is triggered for $ \mathcal{N}_t^\text{FRAG} $. Based on the result of the defragmentation target, $ \mathcal{N}_t^\text{DEFRAG} $ (obtained via \lr{} or \minRWA{}), lightpaths are rerouted. At the end of stage $ t $, the lightpaths associated with the connections leaving the network are dismantled, the uncertainty in the traffic is revealed, and new connection requests are granted using spare resources in the form of available lightpaths via the deterministic \maxRWA{} problem. This process leads to a new fragmented provisioning. 
Note that for simulation purposes, our planning horizon is finite, and defragmentation event times are fixed and given.

By performing an SAA analysis, we choose $ |\xisetSAA|=1000 $ scenarios for \abilene{} and  $ |\xisetSAA|=500 $ for \cost{}, each resulting in worst-case optimality gaps of less than 5\%.   
\revacc{{\App} (Section \ref{appendix:saa-results})}
provides the full SAA results for \lr{}. \rev{Other parameters, namely $ \arrivals $ and $ |\waveSet| $, are selected  such that obtained instances are non-trivial (not too ``free'' or ``saturated'' at the time of the defragmentation), while their size remains small and computationally affordable in our simulation framework. Note that, in an underused network, regardless of the target of the defragmentation procedure (provisionings of \lr{} or \minRWA{}), we can continue granting the incoming requests as resources exist on the network. Furthermore, since rerouting the lightpaths is a costly procedure, defragmentation is usually reserved for the situations where some performance measure has dropped significantly, which is not the case in an under-utilized network. Therefore, in such instances, we are facing a ``trivial'' case where there is virtually no value of performing a defragmentation. On the other hand, a saturated network, one close to its maximum capacity, will not benefit from lightpath rerouting either and needs to seek other enhancements such as network upgrades.} We consider two initial states (denoted by ``Init.'') of 90 and 100 existing requests for \abilene{}, when the network has 6 and 7 wavelengths, respectively. For \cost{}, we use 250 and 260 existing requests over 10 wavelengths in instance generation. Finally, for both networks, we have $ \arrivals=10 $ and $ \drops \in \{15,20\} $.

\begin{figure}[b]
	\centering
	\small
\begin{tikzpicture}

\definecolor{color0}{rgb}{0.12156862745098,0.466666666666667,0.705882352941177}

\begin{groupplot}[width=0.28\textwidth,height=0.3\textwidth,group style={group size=4 by 1}]
\nextgroupplot[
tick align=outside,
tick pos=left,
title={\(\displaystyle \drops=\)15 - Init=90},
x grid style={white!69.0196078431373!black},
xlabel={Stage},
xmajorgrids,
xmin=0.55, xmax=10.45,
xtick style={color=black},
y grid style={white!69.0196078431373!black},
ylabel={Relative granted requests (\%)},
ymajorgrids,
ymin=-15, ymax=40,
ytick style={color=black}
]
\addplot [semithick, color0, mark=*, mark size=2, mark options={solid}]
table {%
1 36.8421052631579
2 12.9629629629629
3 -6.46387832699619
4 5.42986425339367
5 8.00000000000001
6 8.37696335078534
7 -0.990099009900991
8 -7.10659898477157
9 -0.956937799043063
10 1.98412698412698
};
\addplot [semithick, black]
table {%
0.550000000000001 0
10.45 0
};

\nextgroupplot[
tick align=outside,
tick pos=left,
title={\(\displaystyle \drops=\)20 - Init=90},
x grid style={white!69.0196078431373!black},
xlabel={Stage},
xmajorgrids,
xmin=0.55, xmax=10.45,
xtick style={color=black},
y grid style={white!69.0196078431373!black},
ymajorgrids,
ymin=-15, ymax=40,
ytick style={color=black}
]
\addplot [semithick, color0, mark=*, mark size=2, mark options={solid}]
table {%
1 32.5657894736842
2 10.7806691449814
3 -6.79245283018867
4 -1.8348623853211
5 12.2549019607843
6 1.04712041884817
7 0
8 -4.65116279069768
9 -8.53658536585365
10 -4.14201183431952
};
\addplot [semithick, black]
table {%
0.550000000000004 0
10.45 0
};

\nextgroupplot[
tick align=outside,
tick pos=left,
title={\(\displaystyle \drops=\)15 - Init=100},
x grid style={white!69.0196078431373!black},
xlabel={Stage},
xmajorgrids,
xmin=0.55, xmax=10.45,
xtick style={color=black},
y grid style={white!69.0196078431373!black},
ymajorgrids,
ymin=-15, ymax=40,
ytick style={color=black}
]
\addplot [semithick, color0, mark=*, mark size=2, mark options={solid}]
table {%
1 15.3846153846154
2 20.0716845878136
3 7.81893004115226
4 11.1111111111111
5 7.7922077922078
6 5.52763819095478
7 7.0754716981132
8 -1.86046511627907
9 3.41463414634148
10 -1.52671755725191
};
\addplot [semithick, black]
table {%
0.550000000000001 0
10.45 0
};

\nextgroupplot[
tick align=outside,
tick pos=left,
title={\(\displaystyle \drops=\)20 - Init=100},
x grid style={white!69.0196078431373!black},
xlabel={Stage},
xmajorgrids,
xmin=0.55, xmax=10.45,
xtick style={color=black},
y grid style={white!69.0196078431373!black},
ymajorgrids,
ymin=-15, ymax=40,
ytick style={color=black}
]
\addplot [semithick, color0, mark=*, mark size=2, mark options={solid}]
table {%
1 22.4615384615385
2 13.6200716845878
3 7.46887966804978
4 13.6563876651982
5 3.04347826086958
6 -5.16431924882629
7 3.68421052631579
8 -1.60427807486631
9 8.6092715231788
10 -7.18232044198896
};
\addplot [semithick, black]
table {%
0.550000000000004 0
10.45 0
};
\end{groupplot}

\end{tikzpicture}%
	\caption{Relative performance of \lr{} for defragmentation of \abilene{} compared to \minRWA{}}
	\label{fig:defrag-abilene-add}
\end{figure}
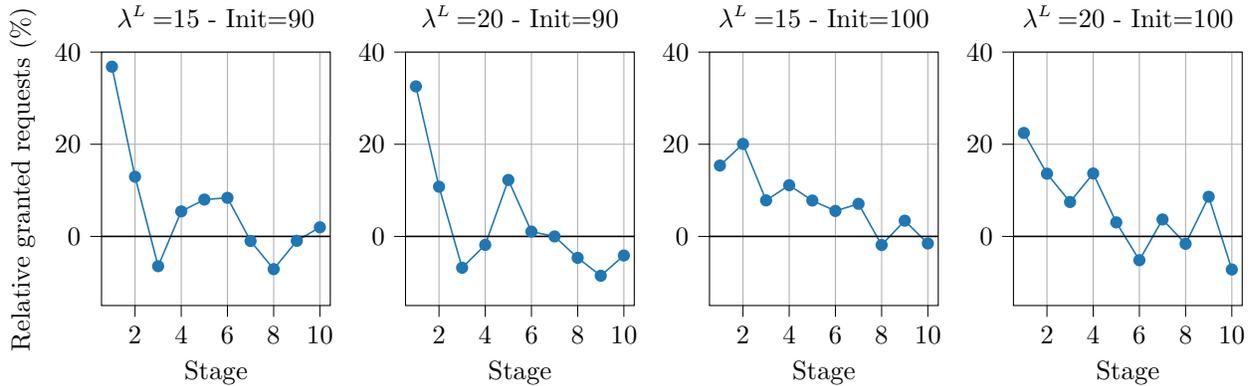
\begin{figure}[htbp]
	\centering
	\small
\begin{tikzpicture}

\definecolor{color0}{rgb}{0.12156862745098,0.466666666666667,0.705882352941177}

\begin{groupplot}[width=0.28\textwidth,height=0.3\textwidth,group style={group size=4 by 1}]
\nextgroupplot[
tick align=outside,
tick pos=left,
title={\(\displaystyle \drops=\)15 - Init=250},
x grid style={white!69.0196078431373!black},
xlabel={Stage},
xmajorgrids,
xmin=0.55, xmax=10.45,
xtick style={color=black},
y grid style={white!69.0196078431373!black},
ylabel={Relative granted requests (\%)},
ymajorgrids,
ymin=-15, ymax=35,
ytick style={color=black}
]
\addplot [semithick, color0, mark=*, mark size=2, mark options={solid}]
table {%
1 29.8550724637681
2 3.02197802197802
3 0.52910052910053
4 0
5 2.86624203821656
6 3.62318840579711
7 3.59281437125749
8 1.13960113960114
9 4.36781609195403
10 4.01337792642139
};
\addplot [semithick, black]
table {%
0.550000000000001 0
10.45 0
};

\nextgroupplot[
tick align=outside,
tick pos=left,
title={\(\displaystyle \drops=\)20 - Init=250},
x grid style={white!69.0196078431373!black},
xlabel={Stage},
xmajorgrids,
xmin=0.55, xmax=10.45,
xtick style={color=black},
y grid style={white!69.0196078431373!black},
ymajorgrids,
ymin=-15, ymax=35,
ytick style={color=black}
]
\addplot [semithick, color0, mark=*, mark size=2, mark options={solid}]
table {%
1 23.1884057971014
2 4.6831955922865
3 2.91777188328912
4 -1.69014084507041
5 2.91262135922331
6 4.15094339622641
7 4.18118466898955
8 7.14285714285714
9 4.93827160493826
10 5.76923076923077
};
\addplot [semithick, black]
table {%
0.550000000000004 0
10.45 0
};

\nextgroupplot[
tick align=outside,
tick pos=left,
title={\(\displaystyle \drops=\)15 - Init=260},
x grid style={white!69.0196078431373!black},
xlabel={Stage},
xmajorgrids,
xmin=0.55, xmax=10.45,
xtick style={color=black},
y grid style={white!69.0196078431373!black},
ymajorgrids,
ymin=-15, ymax=35,
ytick style={color=black}
]
\addplot [semithick, color0, mark=*, mark size=2, mark options={solid}]
table {%
1 25.3462603878116
2 9.78366160800776
3 10.8855224113697
4 -6.4586244371992
5 10.753502461189
6 -0.231347599768647
7 5.26315789473685
8 5.93220338983051
9 4.15382719092816
10 -1.89289012003694
};
\addplot [semithick, black]
table {%
0.550000000000001 0
10.45 0
};

\nextgroupplot[
tick align=outside,
tick pos=left,
title={\(\displaystyle \drops=\)20 - Init=260},
x grid style={white!69.0196078431373!black},
xlabel={Stage},
xmajorgrids,
xmin=0.55, xmax=10.45,
xtick style={color=black},
y grid style={white!69.0196078431373!black},
ymajorgrids,
ymin=-15, ymax=35,
ytick style={color=black}
]
\addplot [semithick, color0, mark=*, mark size=2, mark options={solid}]
table {%
1 26.5664160401003
2 21.6275984077842
3 8.35058661145617
4 -1.81268882175227
5 10.1964949548593
6 -10.3747963063552
7 -1.65631469979295
8 5.66073819085866
9 0.223850223850234
10 4.22841335704405
};
\addplot [semithick, black]
table {%
0.550000000000004 0
10.45 0
};
\end{groupplot}

\end{tikzpicture}%
	\caption{Relative performance of \lr{} for defragmentation of \cost{} compared to \minRWA{}}
	\label{fig:defrag-cost239-add}
\end{figure}
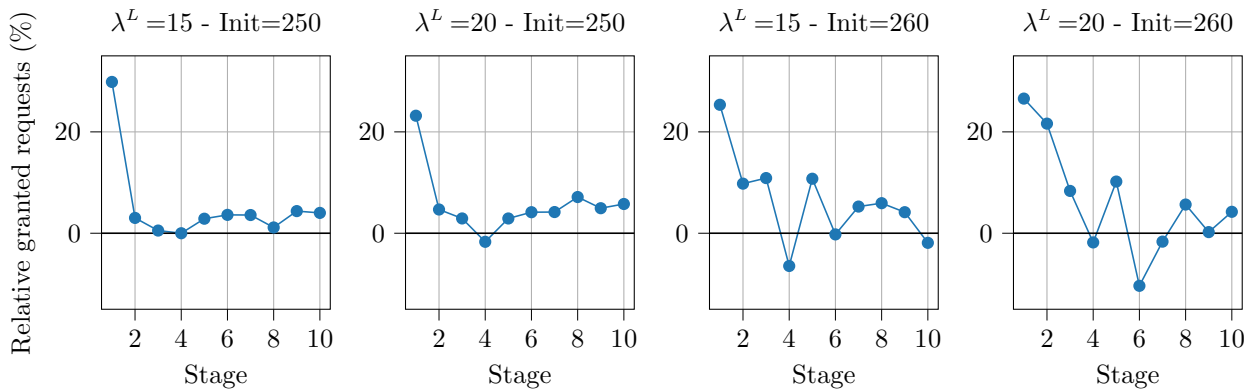

Figures \ref{fig:defrag-abilene-add} and \ref{fig:defrag-cost239-add} show that for both networks, irrespective of the requests' service holding time rate, defragmentation by \lr{} grants more requests over the planning horizon, compared to  \minRWA{}. The improvements are  up to 36\% for \abilene{} and 29\% for \cost{}, each having 90 and 250 initial requests, respectively, with  $ \lambda^L = 15 $. In order to achieve this, \lr{} uses slightly more resources: up to 4\% and 2.8\% more for \abilene{} and \cost{}, respectively (spectrum usage results are provided in 
\revacc{{\App}, Section \ref{appendix:simu-results}).}
This shows that merely minimizing the number of links in defragmentation without considering the future traffic does not necessarily give the best GoS in longer term.

Since the simulation results for \lr{} is not as compelling as \srwa{} simulations, we further investigate the importance of considering the stochasticity of the traffic in our rerouting decisions. This is often achieved by computing the Value of Stochastic Solution (VSS). Let the \emph{mean-value scenario} be the scenario built by using the expected value of random parameters, and  a \emph{mean-value solution} be an optimal solution to the deterministic model employing only the mean-value scenario. VSS is defined as the difference between the optimal value of the stochastic model and the expected long-run objective if the mean-value solution is implemented. Due to the specific discrete nature of the probability distribution of request arrivals, the mean-value scenario is not meaningful in our case. Therefore we use the notion of Expected VSS (EVSS) as defined by \cite{mousavi276stochastic} for the \IPIP{} model. To compute EVSS, we first solve the individual deterministic models associated with every scenario in our sample $ \xisetSAA $, then evaluate the first-stage solutions of each scenario-specific problem using the scenarios in $ \xisetSAA $. EVSS is the difference between the optimal value of the stochastic model and the average of these $ |\xisetSAA| $ scenario objective values.  
Our experiments show that, for the instances where the extensive form of \IPIP{} model is solvable by CPLEX, there is no gap between the optimal values of \IPIP{} and \IPLP{} models (while their solutions can be different). 
Therefore we consider the optimal objective value of the \IPLP{} model as an approximation of  the \IPIP{} model. The EVSS results are presented in Table \ref{tab:evss}, where EVSS numbers (in terms of the number of requests) are reported for different levels of initial requests on the network.
\begin{table}[htbp]
  \centering
  \caption{Expected value of stochastic solution for \lr{}}
\small
\begin{tabular}{ccccccc}
\toprule
\multirow{2}{*}{Network} & \multicolumn{3}{c}{$\arrivals=10$} & \multicolumn{3}{c}{$\arrivals=20$} \\
\cmidrule(lr){2-4}\cmidrule(lr){5-7}       & Init.  & EVSS  & $\sigma^\text{det}$ & Init.  & EVSS  & $\sigma^\text{det}$ \\
\midrule
\multirow{8}{*}{\abilene{}} & \multicolumn{1}{r}{80} & 0.71  & 0.44  & \multicolumn{1}{r}{80} & 1.35  & 0.65 \\
      & \multicolumn{1}{r}{90} & 0.77  & 0.65  & \multicolumn{1}{r}{90} & 1.60  & 1.1 \\
      & \multicolumn{1}{r}{100} & 1.48  & 0.64  & \multicolumn{1}{r}{100} & 2.52  & 0.97 \\
      & \multicolumn{1}{r}{110} & 2.92  & 0.67  & \multicolumn{1}{r}{110} & 3.86  & 0.84 \\
      & \multicolumn{1}{r}{120} & 1.40  & 0.38  & \multicolumn{1}{r}{120} & 2.03  & 0.57 \\
      & \multicolumn{1}{r}{130} & 1.19  & 0.24  & \multicolumn{1}{r}{130} & 1.52  & 0.45 \\
      & \multicolumn{1}{r}{140} & 0.21  & 0.09  & \multicolumn{1}{r}{140} & 0.40  & 0.19 \\
      & \multicolumn{1}{r}{150} & 0.06  & 0.06  & \multicolumn{1}{r}{150} & 0.11  & 0.11 \\
\midrule
\multirow{6}{*}{\cost{}} & 200   & 1.57  & 0.60  & 200   & 2.85  & 1.06 \\
      & 210   & 1.48  & 0.67  & 210   & 2.93  & 1.24 \\
      & 220   & 2.66  & 0.84  & 220   & 4.72  & 1.37 \\
      & 230   & 3.06  & 0.77  & 230   & 5.21  & 1.46 \\
      & 240   & 4.04  & 0.79  & 240   & 6.66  & 1.48 \\
      & 250   & 3.04  & 0.46  & 250   & 4.43  & 0.74 \\
\bottomrule
\end{tabular}%
  \label{tab:evss}%
\end{table}%
$ \sigma^\text{det} $ is the standard deviation of the evaluation over $ |\xisetSAA| $ of deterministic subproblems associated with each scenario in $ \xisetSAA $. The results show that the stochastic model attains significant gains of up to 3 and 4 more requests out of an average of 10 incoming arrivals in \abilene{} and \cost{}, respectively, showcasing the value of the stochastic model over a deterministic approach.  

Considerable benefit from solving the stochastic model suggested by EVSS analysis is not reflected in the moderate performance of \lr{} in simulation results. One reason could be the  \IPLP{} approximation. Although in all the tested instances that we were able to solve \IPIP{} models (either using the extensive form or the integer L-shaped method), there is no gap between \IPIP{} and \IPLP{} optimal values, but their solutions are quite different. The first-stage solutions from the \IPLP{} model are not necessarily optimal for the \IPIP{} model. Another reason can be the intrinsically multistage nature of the problem that is approximated by \lr{}, a two-stage model. \lr{} looks one stage ahead into the future in its rerouting decisions, while the traffic arrives and leaves the network over the entire planning horizon. Figures \ref{fig:defrag-abilene-add} and \ref{fig:defrag-cost239-add} support this claim by showing that the most significant improvements in terms of granted requests are obtained at the early stages. Finally, note that, as a major source of network fragmentation, drops in the simulation framework are not modeled by \lr{}, which could lead to losing more benefits. 

\section{Conclusion}
\label{sec:conlusion}
Stochastic nature of the traffic in WDM networks is a contributing factor in sub-optimal usage of the available bandwidth.  To tackle this issue, in this paper we have introduced \srwa{}, the first 2SP model for \maxRWA{} as the provisioning problem of WDM networks, with the purpose of making more informed lightpath assignment decisions. Furthermore, for the costly operation of network defragmentation that aims to improve the spectrum usage in WDM networks, we have proposed \lr{}, a 2SP lightpath rerouting model as the provisioning target of the procedure. 

In order to solve the 2SP models, we have designed a decomposition framework based on two different relaxations that make use of our newly developed problem-specific family of cuts. Our numerical results show that these cuts are quite effective in reducing the number of Benders optimality cuts needed for the algorithm to converge, in turn improving the scalability of the solution framework. 
Through a set of comprehensive numerical experiments over six standard long-haul networks in the literature, we have showed that \srwa{}  considerably improves the overall GoS in all instances compared to the results from deterministic \maxRWA{}. At the same time, \srwa{} reduces the spectrum usage despite having a higher number of connection requests on the network, leaving more room for future traffic.  
As for the  \lr{} problem, our numerical experiments are promising, showing moderate improvements in the number of granted requests and overall GoS. 

Our numerical study indicates the merit of considering traffic uncertainty in the decision making involved in provisioning and defragmentation of WDM networks. This leads to several questions that open new areas of research. First, our 2SP models  look one step ahead into the future traffic, while the multistage nature of the provisioning and defragmentation decisions suggests that a multistage model might better capture the  process.  Multistage problems involving integer decision variables are notoriously difficult to solve, therefore the practicality of them for \srwa{} and \lr{} problems remains open. Second, although our designed solution method is able to demonstrate the value of using a stochastic model for several networks, it is not able to solve very large instances where  network size, arrival of the connection requests or the number of available wavelengths is considerably larger than the studied instances. In addition, as demonstrated by our evaluation of \lr{} problem, solving the exact \IPIP{} model is also of interest. Therefore, new methodological contributions are needed. \rev{For the \lr{} problem, another interesting research direction is to integrate the lightpath rerouting decisions into the provisioning problem such that the legacy provisioning can be migrated to the new one without (or with minimum number of) disruptions.} Finally,  significant savings from \srwa{} problem for WDM networks unfold new questions about the impact that ignoring the future traffic can have on other practical decisions concerning WDM and more generally other types of optical networks. For instance,  an interesting question is whether a stochastic model can determine the best time for triggering a defragmentation process by taking the future traffic into account. As another example,  flex-grid networks, the popular emerging  transmission systems, face the same traffic uncertainty as WDM networks. Thus, their provisioning, a much more difficult problem than the RWA, might also benefit from stochastic models. 
 
\newpage
\ACKNOWLEDGMENT{%
This work was supported by Natural Sciences and Engineering Research Council of Canada [Grant RGPIN- 2018-04984].
Computations were performed on the Niagara supercomputer at the SciNet HPC Consortium. SciNet is funded by: the Canada Foundation for Innovation; the Government of Ontario; Ontario Research Fund - Research Excellence; and the University of Toronto.
}

\bibliographystyle{informs2014} 
\bibliography{Bibliography.bib} 

\begin{appendices}

\section{Extensive Form for the \srwa{}  Problem} \label{appendix:EXT-form}
Denote by $ \FeasR{x,\xival,\SD} $ the feasible set of $\Q(x,\xival,\SD)$, defined in the paper by \eqref{eq: conflict_link_demands_scenarios}-\eqref{eq: granted_scenarios_bounds}. We can reformulate the SAA problem \eqref{eq: saa-SmaxRWA} as a mixed-integer programming model:
\BSE
\label{eq: saa-EXF-SmaxRWA}
\begin{alignat}{5}
	\max \ \ & \sum_{(s,d) \in \SDcurr{n}}\sum_{\wavelength \in \waveSet}\sum_{\stackrel{\edge\in}{\nodefunc^+(s)\cap \edgeSetw}}\link + \frac{1}{|\xisetSAA|} \sum_{\xival \in \xisetSAA} \eta_\xival \label{eq: OF_EXF_SAA}\\
	\text{s.t.}\ \ & x\in \Feas{\SDcurr{n}} \\
	& \eta_\xival \leq \sum_{(s,d) \in \SDsce} \numgrantedsce(\xival) &\qquad& \xival \in \xisetSAA\\
	& (y(\xival),z(\xival)) \in \FeasR{x,\xival,\SDcurr{n}} &\qquad& \xival \in \xisetSAA\\
	& \eta_\xival \geq 0 && \xival \in \xisetSAA, 
\end{alignat}
\ESE
i.e., by creating copies of the second-stage decision variables and constraints for each scenario in the sample. This is  the \emph{extensive form} of the SAA problem.

\section{SAA Analysis}\label{appendix:saa-results}
In an SAA analysis, we decide on the sample size to be used in the SAA problem, in order to get statistically valid lower and upper bound on the optimal value of the 2SP model. Consider a  maximization problem. In an SAA analysis with $ n $ repetitions, for each scenario level $ |\xisetSAA| $, we generate $ n $ samples of size $ |\xisetSAA| $ and solve  the SAA problem for each sample. The average of the optimal values of these SAA problems provides a statistically valid upper bound. We then evaluate the solutions of each SAA problem with a much larger evaluation sample and pick the best one, which gives us a statistical lower bound. This process is repeated for different scenario levels, until a predefined condition is met.  For the details of the steps involved in this analysis, interested reader can refer to \cite{shapiro2014lectures}. In the following, for both  \srwa{} and \lr{}, we have performed 
an SAA analysis considering the two smallest networks, \abilene{} and \cost{}, with 30 repetitions and an evaluation sample of 5000 scenarios.  

For \srwa{}, we have considered two scenario levels of $ |\xisetSAA|\in\{100,500\} $, as the larger sample sizes failed to converge in the given time limit. A contributing factor in this analysis can be the state of the network when \srwa{} is solved, i.e., which wavelinks are already serving connections on the network. For generating such initial states, we have considered the simulation framework given in Section \ref{sec:simul-SmaxRWA} of the paper, and following one sample path (with \maxRWA{} used to grant new arrivals), we have taken three snapshots of the state of the network at stages 13, 26 and 52, to be used as the initial state. Then, for each combination of network, stage and scenario level, we have solved the instance we have performed an SAA analysis with 30 repetitions.  Results are summarized in Table \ref{tab:app-saa-SmaxRWA}, where the mean and width of the confidence intervals for upper bound (UB) and lower bound (LB) results are given, as well as the final worst-case optimality gap.  Based on these results,  for \abilene{} and \cost{} instances, we have chosen sample size of $ |\xisetSAA|=500 $  where the optimality gap is less than or equal to 3\%. Note that, our main goal is to retrieve feasible solutions to be used in the simulation framework, hence the small LB widths are adequate. 
\begin{table}[htbp]
	\centering
	\caption{\srwa{} SAA results}
	\small
	\begin{tabular}{ccrrrrrr}
		\toprule
		\multirow{2}{*}{Network} & \multirow{2}{*}{Stage} & \multicolumn{1}{c}{\multirow{2}{*}{$|\xisetSAA|$}} & \multicolumn{2}{c}{UB} & \multicolumn{2}{c}{LB} & \multicolumn{1}{c}{\multirow{2}{*}{gap (\%)}} \\
		\cmidrule(lr){4-5}\cmidrule(lr){6-7}          &       &       & \multicolumn{1}{c}{Mean} & \multicolumn{1}{c}{Width} & \multicolumn{1}{c}{Mean} & \multicolumn{1}{l}{Width} &  \\
		\midrule
		\multirow{6}{*}{\abilene{}} & \multirow{2}{*}{13} & 100   & 19.63 & 0.07  & 19.10 & 0.11  & 4.36 \\
		&       & 500   & 19.70 & 0.04  & 19.23 & 0.02  & 2.18 \\
		\cmidrule{2-8}          & \multirow{2}{*}{26} & 100   & 21.88 & 0.07  & 21.30 & 0.16  & 5.08 \\
		&       & 500   & 21.89 & 0.05  & 21.47 & 0.18  & 2.35 \\
		\cmidrule{2-8}          & \multirow{2}{*}{52} & 100   & 21.44 & 0.07  & 21.10 & 0.11  & 4.17 \\
		&       & 500   & 21.45 & 0.07  & 21.19 & 0.07  & 2.06 \\
		\midrule
		\multirow{6}{*}{\cost{}} & \multirow{2}{*}{13} & 100   & 21.02 & 0.09  & 20.40 & 0.18  & 3.71 \\
		&       & 500   & 20.36 & 0.04  & 20.09 & 0.13  & 2.73 \\
		\cmidrule{2-8}          & \multirow{2}{*}{26} & 100   & 19.03 & 0.09  & 18.37 & 0.17  & 3.83 \\
		&       & 500   & 19.12 & 0.03  & 18.72 & 0.02  & 3.09 \\
		\cmidrule{2-8}          & \multirow{2}{*}{52} & 100   & 21.96 & 0.09  & 21.33 & 0.17  & 2.44 \\
		&       & 500   & 22.03 & 0.07  & 21.70 & 0.04  & 1.93 \\
		\bottomrule
	\end{tabular}%
	\label{tab:app-saa-SmaxRWA}%
\end{table}%

For each network considered in  the \lr{} problem, two levels of initial connection requests on the network (denoted by ``Init.'') are considered. Requests are then randomly generated while making sure that they can result in a feasible provisioning. Table \ref{tab:app-saa-SmaxLR} presents the SAA results for three scenario levels of $ |\xisetSAA|\in\{10,500,1000\} $. Based on these experiments, for \abilene{} and \cost{} we choose $ \xisetSAA = 1000 $ and $ \xisetSAA = 500 $, respectively, where at this level the worst-case optimality gap is less than 5\% for instances of both networks.

\renewcommand{\aboverulesep}{0pt}
\renewcommand{\belowrulesep}{0pt}
\begin{table}[htbp]
	\centering
	\caption{\lr{} SAA results}
	\small
	\begin{tabular}{ccrrrrrr}
		\toprule
		\multirow{2}{*}{Network} & \multirow{2}{*}{Init.} & \multicolumn{1}{c}{\multirow{2}{*}{$|\xisetSAA|$}} & \multicolumn{2}{c}{UB CI} & \multicolumn{2}{c}{LB CI} & \multicolumn{1}{c}{\multirow{2}{*}{gap (\%)}} \\
		\cmidrule(lr){4-5}\cmidrule(lr){6-7}      &       &       & \multicolumn{1}{c}{Mean} & \multicolumn{1}{c}{Width} & \multicolumn{1}{c}{Mean} & \multicolumn{1}{l}{Width} &  \\
		\midrule
		\multirow{6}{*}{\abilene{}} & \multirow{3}{*}{100} & 10    & 9.97  & 0.31  & 9.04  & 0.05  & 14.36 \\
		&       & 500   & 9.98  & 0.04  & 9.71  & 0.08  & 4.10 \\
		&       & 1000  & 9.97  & 0.03  & 9.89  & 0.08  & 2.03 \\
		\cmidrule{2-8}      & \multirow{3}{*}{120} & 10    & 4.59  & 0.24  & 4.09  & 0.03  & 19.00 \\
		&       & 500   & 4.69  & 0.03  & 4.48  & 0.04  & 6.18 \\
		&       & 1000  & 4.68  & 0.02  & 4.53  & 0.04  & 4.79 \\
		\midrule
		\multirow{6}{*}{\cost{}} & \multirow{3}{*}{200} & 10    & 9.83  & 0.39  & 9.62  & 0.08  & 7.24 \\
		&       & 500   & 9.97  & 0.05  & 9.85  & 0.08  & 2.57 \\
		&       & 1000  & 9.99  & 0.03  & 9.99  & 0.09  & 1.24 \\
		\cmidrule{2-8}      & \multirow{3}{*}{220} & 10    & 9.84  & 0.35  & 9.52  & 0.08  & 7.97 \\
		&       & 500   & 10.03 & 0.04  & 9.77  & 0.08  & 3.87 \\
		&       & 1000  & 10.02 & 0.03  & 10.01 & 0.09  & 1.30 \\
		\bottomrule
	\end{tabular}%
	\label{tab:app-saa-SmaxLR}%
\end{table}%

\section{Simulation Results}\label{appendix:simu-results}
In this section, we present the simulation results that are removed from the main text for the sake of brevity. 

\subsection{Performance}
\begin{figure}[htbp]
	\centering
	\small
	\input{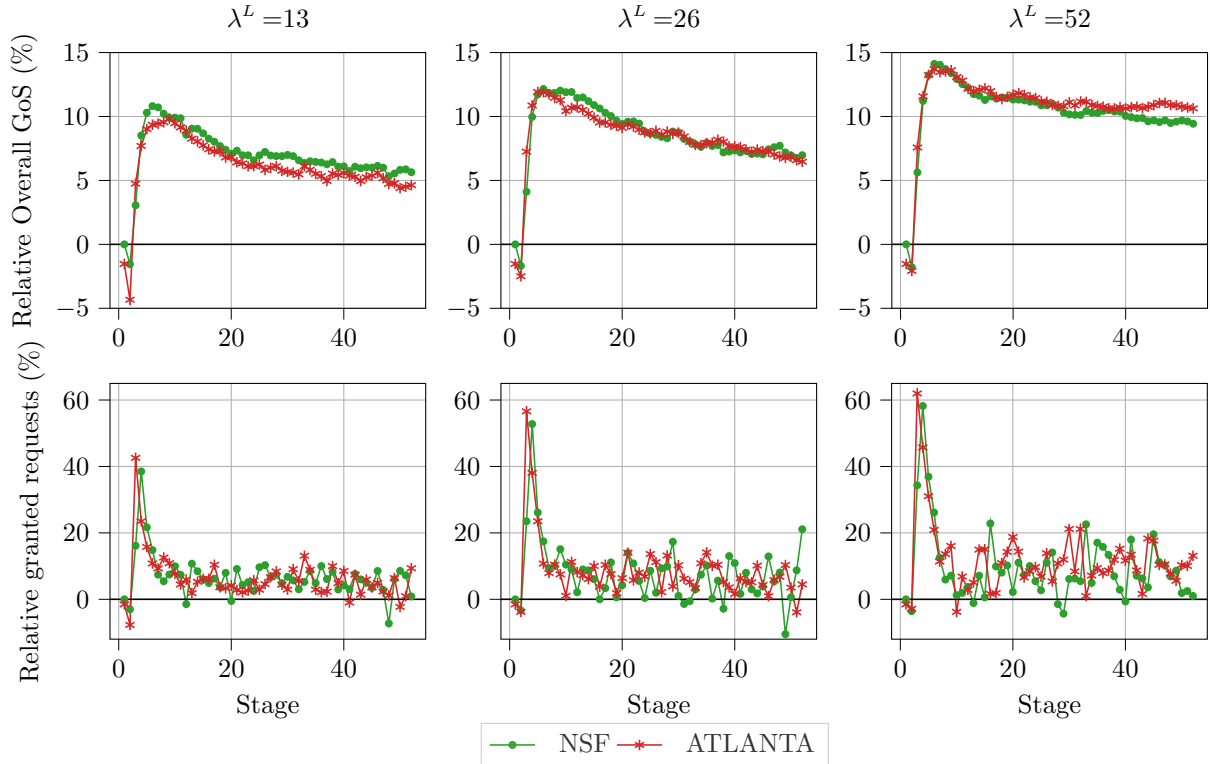}%
	\caption{Relative performance of \srwa{} for medium-sized networks over different simulation horizons.}
	\label{fig:app-improved-nsfatlantagos-add}
\end{figure}

Figure \ref{fig:app-improved-nsfatlantagos-add} illustrates the relative improvements obtained from employing \srwa{} over the \maxRWA{} problem for medium-sized networks \nsf{} and \atlanta{}. Together with Figures \ref{fig:improved-cost239abilenegos-add} and \ref{fig:improved-usabrazilgos-add} given in the paper, they provide a general picture of the merit of incorporating stochasticity in our decision-making in terms of granting requests over different network.

\subsection{Spectrum Usage}
In addition to the results presented for \abilene{}, \cost{} and \brazil{} in Figure \ref{fig:spec-usage} of the paper, we have analyzed the spectrum usage of \nsf{}, \atlanta{} and \usa{} and the results are depicted in  Figure \ref{fig:app-spec-usage}. The results are consistent with previously mentioned network.

\begin{figure}
	\centering
	\hspace*{-4mm}
	\subfigure[\nsf{}]{
		\centering
		\small
		\input{./Figures/improved-NSFlinkusage.tex}
		\label{fig:app-spectrum-usage-NSF}
	}
	
	\hspace*{-4mm}
	\subfigure[\atlanta{}]{
		\centering
		\small
		\input{./Figures/improved-atlantalinkusage.tex}
		\label{fig:app-spectrum-usage-atlanta}
	}
	
	\hspace*{-4mm}
	\subfigure[\usa{}]{
		\centering
		\small
		\input{./Figures/improved-USAlinkusage.tex}
		\label{fig:app-spectrum-usage-USA}
	}

	\caption{Relative bandwidth spectrum usage}
	\label{fig:app-spec-usage}
\end{figure}

\begin{figure}
	\centering
	\hspace*{-4mm}
	\subfigure[\abilene{}]{
		\centering
		\small
\begin{tikzpicture}

\definecolor{color0}{rgb}{1,0.498039215686275,0.0549019607843137}

\begin{axis}[
width=0.7\linewidth,
height=0.3\linewidth,
legend cell align={left},
legend columns=2,
legend style={fill opacity=0.8, draw opacity=1, text opacity=1, at={(0.5,1.4)}, anchor=north, draw=white!80!black},
tick align=outside,
tick pos=left,
x grid style={white!69.0196078431373!black},
xlabel={Stage},
xmin=0.41, xmax=11.19,
xtick style={color=black},
xtick={1.22,3.22,5.22,7.22,9.22},
xticklabels={1,3,5,7,9},
y grid style={white!69.0196078431373!black},
ylabel={Relative spectrum usage (\%)},
ymin=-5, ymax=0,
ytick style={color=black}
]
\draw[draw=black] (axis cs:0.9,0) rectangle (axis cs:1.1,-4.04390525707685);
\addlegendimage{ybar,ybar legend,draw=black};
\addlegendentry{$\drops=$15 - Init=90$\quad$}

\draw[draw=black] (axis cs:1.9,0) rectangle (axis cs:2.1,-4.00701754385965);
\draw[draw=black] (axis cs:2.9,0) rectangle (axis cs:3.1,-3.42442298472707);
\draw[draw=black] (axis cs:3.9,0) rectangle (axis cs:4.1,-3.42945694650317);
\draw[draw=black] (axis cs:4.9,0) rectangle (axis cs:5.1,-3.48914635770408);
\draw[draw=black] (axis cs:5.9,0) rectangle (axis cs:6.1,-3.52289220213747);
\draw[draw=black] (axis cs:6.9,0) rectangle (axis cs:7.1,-3.38649340421343);
\draw[draw=black] (axis cs:7.9,0) rectangle (axis cs:8.1,-3.09685054901703);
\draw[draw=black] (axis cs:8.9,0) rectangle (axis cs:9.1,-3.25692012948405);
\draw[draw=black] (axis cs:9.9,0) rectangle (axis cs:10.1,-3.44088900789933);
\draw[draw=black,fill=color0,fill opacity=0.3] (axis cs:1.1,0) rectangle (axis cs:1.3,-3.4734257654535);
\addlegendimage{ybar,ybar legend,draw=black,fill=color0,fill opacity=0.3};
\addlegendentry{$\drops=$20 - Init=90$\quad$}

\draw[draw=black,fill=color0,fill opacity=0.3] (axis cs:2.1,0) rectangle (axis cs:2.3,-3.61902090054706);
\draw[draw=black,fill=color0,fill opacity=0.3] (axis cs:3.1,0) rectangle (axis cs:3.3,-3.14060896339378);
\draw[draw=black,fill=color0,fill opacity=0.3] (axis cs:4.1,0) rectangle (axis cs:4.3,-3.05841231431067);
\draw[draw=black,fill=color0,fill opacity=0.3] (axis cs:5.1,0) rectangle (axis cs:5.3,-3.13348790990395);
\draw[draw=black,fill=color0,fill opacity=0.3] (axis cs:6.1,0) rectangle (axis cs:6.3,-2.83432611114747);
\draw[draw=black,fill=color0,fill opacity=0.3] (axis cs:7.1,0) rectangle (axis cs:7.3,-2.88886003635421);
\draw[draw=black,fill=color0,fill opacity=0.3] (axis cs:8.1,0) rectangle (axis cs:8.3,-2.64563887994851);
\draw[draw=black,fill=color0,fill opacity=0.3] (axis cs:9.1,0) rectangle (axis cs:9.3,-2.58167841127482);
\draw[draw=black,fill=color0,fill opacity=0.3] (axis cs:10.1,0) rectangle (axis cs:10.3,-2.52946798343421);
\draw[draw=black,fill=color0,fill opacity=0.7] (axis cs:1.3,0) rectangle (axis cs:1.5,-2.20791721814462);
\addlegendimage{ybar,ybar legend,draw=black,fill=color0,fill opacity=0.7};
\addlegendentry{$\drops=$15 - Init=100$\quad$}

\draw[draw=black,fill=color0,fill opacity=0.7] (axis cs:2.3,0) rectangle (axis cs:2.5,-2.80555228796611);
\draw[draw=black,fill=color0,fill opacity=0.7] (axis cs:3.3,0) rectangle (axis cs:3.5,-3.03503295178633);
\draw[draw=black,fill=color0,fill opacity=0.7] (axis cs:4.3,0) rectangle (axis cs:4.5,-3.08279674055502);
\draw[draw=black,fill=color0,fill opacity=0.7] (axis cs:5.3,0) rectangle (axis cs:5.5,-2.96450469708049);
\draw[draw=black,fill=color0,fill opacity=0.7] (axis cs:6.3,0) rectangle (axis cs:6.5,-3.19404821838115);
\draw[draw=black,fill=color0,fill opacity=0.7] (axis cs:7.3,0) rectangle (axis cs:7.5,-3.11925583467944);
\draw[draw=black,fill=color0,fill opacity=0.7] (axis cs:8.3,0) rectangle (axis cs:8.5,-2.80529889791829);
\draw[draw=black,fill=color0,fill opacity=0.7] (axis cs:9.3,0) rectangle (axis cs:9.5,-2.84817223592733);
\draw[draw=black,fill=color0,fill opacity=0.7] (axis cs:10.3,0) rectangle (axis cs:10.5,-2.61631213741326);
\draw[draw=black,fill=color0] (axis cs:1.5,0) rectangle (axis cs:1.7,-2.79749729274456);
\addlegendimage{ybar,ybar legend,draw=black,fill=color0};
\addlegendentry{$\drops=$20 - Init=100$\quad$}

\draw[draw=black,fill=color0] (axis cs:2.5,0) rectangle (axis cs:2.7,-3.01658238268846);
\draw[draw=black,fill=color0] (axis cs:3.5,0) rectangle (axis cs:3.7,-3.07283544157569);
\draw[draw=black,fill=color0] (axis cs:4.5,0) rectangle (axis cs:4.7,-3.21351382095325);
\draw[draw=black,fill=color0] (axis cs:5.5,0) rectangle (axis cs:5.7,-3.03047309052136);
\draw[draw=black,fill=color0] (axis cs:6.5,0) rectangle (axis cs:6.7,-2.83794290791957);
\draw[draw=black,fill=color0] (axis cs:7.5,0) rectangle (axis cs:7.7,-2.93064507243981);
\draw[draw=black,fill=color0] (axis cs:8.5,0) rectangle (axis cs:8.7,-2.76744694815139);
\draw[draw=black,fill=color0] (axis cs:9.5,0) rectangle (axis cs:9.7,-2.87350054525626);
\draw[draw=black,fill=color0] (axis cs:10.5,0) rectangle (axis cs:10.7,-2.71647779645393);
\addplot [semithick, black, forget plot]
table {%
0.41 0
11.19 0
};
\end{axis}

\end{tikzpicture}
		\label{fig:defrag-spectrum-usage-abilene}
	}
	
	\hspace*{-4mm}
	\subfigure[\cost{}]{
		\centering
		\small
\begin{tikzpicture}

\definecolor{color0}{rgb}{0.12156862745098,0.466666666666667,0.705882352941177}

\begin{axis}[
width=0.7\linewidth,
height=0.3\linewidth,
legend cell align={left},
legend columns=2,
legend style={fill opacity=0.8, draw opacity=1, text opacity=1, at={(0.5,1.4)}, anchor=north, draw=white!80!black},
tick align=outside,
tick pos=left,
x grid style={white!69.0196078431373!black},
xlabel={Stage},
xmin=0.41, xmax=11.19,
xtick style={color=black},
xtick={1.22,3.22,5.22,7.22,9.22},
xticklabels={1,3,5,7,9},
y grid style={white!69.0196078431373!black},
ylabel={Relative spectrum usage (\%)},
ymin=-3, ymax=0,
ytick style={color=black}
]
\draw[draw=black] (axis cs:0.9,0) rectangle (axis cs:1.1,-2.80970625798212);
\addlegendimage{ybar,ybar legend,draw=black};
\addlegendentry{$\drops=$15 - Init=250$\quad$}

\draw[draw=black] (axis cs:1.9,0) rectangle (axis cs:2.1,-2.36404410390581);
\draw[draw=black] (axis cs:2.9,0) rectangle (axis cs:3.1,-2.13615999261004);
\draw[draw=black] (axis cs:3.9,0) rectangle (axis cs:4.1,-2.06187926220879);
\draw[draw=black] (axis cs:4.9,0) rectangle (axis cs:5.1,-1.91586554851026);
\draw[draw=black] (axis cs:5.9,0) rectangle (axis cs:6.1,-1.88187870827483);
\draw[draw=black] (axis cs:6.9,0) rectangle (axis cs:7.1,-1.89703122526964);
\draw[draw=black] (axis cs:7.9,0) rectangle (axis cs:8.1,-1.97009634611889);
\draw[draw=black] (axis cs:8.9,0) rectangle (axis cs:9.1,-2.28441546514831);
\draw[draw=black] (axis cs:9.9,0) rectangle (axis cs:10.1,-2.17726127135257);
\draw[draw=black,fill=color0,fill opacity=0.3] (axis cs:1.1,0) rectangle (axis cs:1.3,-2.58265928764013);
\addlegendimage{ybar,ybar legend,draw=black,fill=color0,fill opacity=0.3};
\addlegendentry{$\drops=$20 - Init=250$\quad$}

\draw[draw=black,fill=color0,fill opacity=0.3] (axis cs:2.1,0) rectangle (axis cs:2.3,-2.34997313648999);
\draw[draw=black,fill=color0,fill opacity=0.3] (axis cs:3.1,0) rectangle (axis cs:3.3,-2.12633328715889);
\draw[draw=black,fill=color0,fill opacity=0.3] (axis cs:4.1,0) rectangle (axis cs:4.3,-1.9479183337525);
\draw[draw=black,fill=color0,fill opacity=0.3] (axis cs:5.1,0) rectangle (axis cs:5.3,-1.746388151777);
\draw[draw=black,fill=color0,fill opacity=0.3] (axis cs:6.1,0) rectangle (axis cs:6.3,-1.77192903254933);
\draw[draw=black,fill=color0,fill opacity=0.3] (axis cs:7.1,0) rectangle (axis cs:7.3,-1.74726359232012);
\draw[draw=black,fill=color0,fill opacity=0.3] (axis cs:8.1,0) rectangle (axis cs:8.3,-1.71885826683654);
\draw[draw=black,fill=color0,fill opacity=0.3] (axis cs:9.1,0) rectangle (axis cs:9.3,-1.81741052490567);
\draw[draw=black,fill=color0,fill opacity=0.3] (axis cs:10.1,0) rectangle (axis cs:10.3,-1.74578047312211);
\draw[draw=black,fill=color0,fill opacity=0.7] (axis cs:1.3,0) rectangle (axis cs:1.5,-2.29199058976471);
\addlegendimage{ybar,ybar legend,draw=black,fill=color0,fill opacity=0.7};
\addlegendentry{$\drops=$15 - Init=260$\quad$}

\draw[draw=black,fill=color0,fill opacity=0.7] (axis cs:2.3,0) rectangle (axis cs:2.5,-2.03745380518358);
\draw[draw=black,fill=color0,fill opacity=0.7] (axis cs:3.3,0) rectangle (axis cs:3.5,-2.24270814299215);
\draw[draw=black,fill=color0,fill opacity=0.7] (axis cs:4.3,0) rectangle (axis cs:4.5,-1.93567424726778);
\draw[draw=black,fill=color0,fill opacity=0.7] (axis cs:5.3,0) rectangle (axis cs:5.5,-2.02516383951487);
\draw[draw=black,fill=color0,fill opacity=0.7] (axis cs:6.3,0) rectangle (axis cs:6.5,-1.88128785797455);
\draw[draw=black,fill=color0,fill opacity=0.7] (axis cs:7.3,0) rectangle (axis cs:7.5,-1.91509609935355);
\draw[draw=black,fill=color0,fill opacity=0.7] (axis cs:8.3,0) rectangle (axis cs:8.5,-1.941417665315);
\draw[draw=black,fill=color0,fill opacity=0.7] (axis cs:9.3,0) rectangle (axis cs:9.5,-2.02105263157895);
\draw[draw=black,fill=color0,fill opacity=0.7] (axis cs:10.3,0) rectangle (axis cs:10.5,-1.9469653479254);
\draw[draw=black,fill=color0] (axis cs:1.5,0) rectangle (axis cs:1.7,-2.34313827044162);
\addlegendimage{ybar,ybar legend,draw=black,fill=color0};
\addlegendentry{$\drops=$20 - Init=260$\quad$}

\draw[draw=black,fill=color0] (axis cs:2.5,0) rectangle (axis cs:2.7,-2.25021697219444);
\draw[draw=black,fill=color0] (axis cs:3.5,0) rectangle (axis cs:3.7,-2.26000938662045);
\draw[draw=black,fill=color0] (axis cs:4.5,0) rectangle (axis cs:4.7,-2.03384759476547);
\draw[draw=black,fill=color0] (axis cs:5.5,0) rectangle (axis cs:5.7,-2.13939086174012);
\draw[draw=black,fill=color0] (axis cs:6.5,0) rectangle (axis cs:6.7,-1.97339134333706);
\draw[draw=black,fill=color0] (axis cs:7.5,0) rectangle (axis cs:7.7,-1.92610377484327);
\draw[draw=black,fill=color0] (axis cs:8.5,0) rectangle (axis cs:8.7,-1.8992574267947);
\draw[draw=black,fill=color0] (axis cs:9.5,0) rectangle (axis cs:9.7,-1.7913453639664);
\draw[draw=black,fill=color0] (axis cs:10.5,0) rectangle (axis cs:10.7,-1.78106946742592);
\addplot [semithick, black, forget plot]
table {%
0.41 -4.44089209850063e-16
11.19 -4.44089209850063e-16
};
\end{axis}

\end{tikzpicture}
		\label{fig:defrag-spectrum-usage-cost239}
	}
	
	\caption{Relative bandwidth spectrum usage after each rerouting}
	\label{fig:app-defrag-spectrum-usage}
\end{figure}
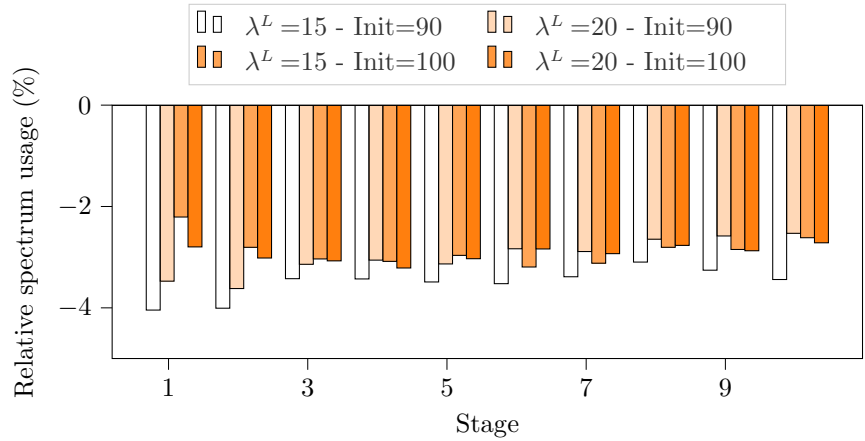
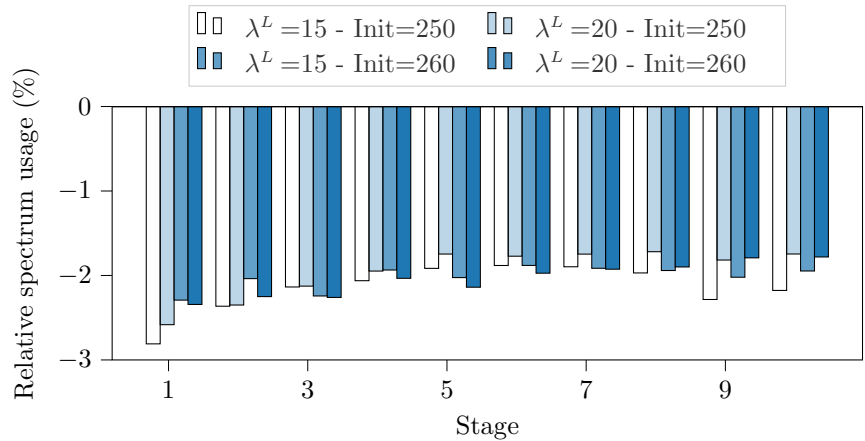

Likewise, Figure \ref{fig:app-defrag-spectrum-usage} depicts the relative spectrum usage by \lr{} compared to the \minRWA{} problem, which is slightly more than that of the latter.

\begin{figure}[htb]
    \centering
    \input{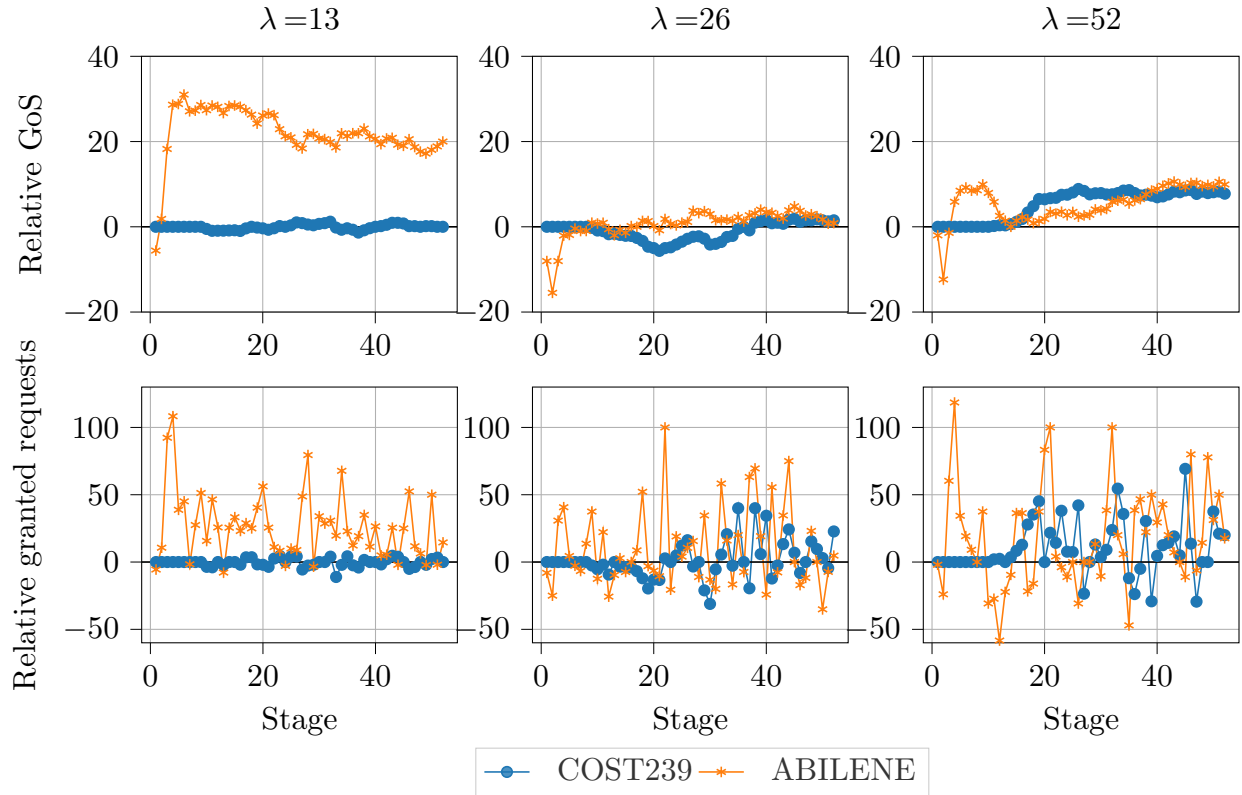}
    \caption{\revminor{Results for 20/80 distribution}}
    \label{fig:2080}
\end{figure}

\revminor{
\subsection{Non-uniform Demand}
We have tested the 20/80 distribution for \cost{} and \abilene{} and the computational results are given in Figure \ref{fig:2080} which illustrates the relative performance of our \srwa{} model compared to the deterministic \maxRWA{}, in terms of the overall GoS and the number of granted requests at each simulated stage. 
}

\revminor{
\indent When the lifetime is short ($\lambda=13$), under the 20/80 distribution in the mesh network  \cost{}, both models perform the same.
The reason is that, compared to a uniform case, the non-uniform distribution yields a less uncertain setting. We expect the traffic to be, with a high probability, confined to a very small subset of the nodes, and when a request is dropped there is a high chance that the incoming request can simply replace it. On the other hand, when the lifetime gets longer, the decisions at early stages become more important for the sake of future requests, as such the stochastic model has potential to perform better. This is observed from the results of $\lambda = 26$ and $\lambda = 52$. More specifically, when $\lambda = 26$, the stochastic model starts to reserve some resources with the hope of achieving a better performance in the future (potentially outside of the planning horizon). When $\lambda = 52$, the decisions matter the most because the resources are released much later in time and stochastic model shows its advantage.

\indent For the \abilene{} network, the stochastic model is still able to outperform the deterministic one, as the resources are rather scarce in comparison to the \cost{} network, which is again a function of the structure of the network rather than the data distribution. We note that for $\lambda = 13$, the GoS performance of both models has stabled after the initial decision-making phase, suggesting again that the new incoming traffic is with a high probability simply replacing the dropped requests.  
}

\end{appendices}

\end{document}